\documentclass[twoside]{irmaems}
\usepackage{amssymb} 
\usepackage{amsmath} 
\usepackage{amscd}
\usepackage{latexsym}
\usepackage{bm}
\usepackage[dvips]{graphicx}

\renewcommand\theequation{\arabic{section}.\arabic{equation}}

\theoremstyle{definition} 

 \newtheorem{definition}{Definition}[section]

\theoremstyle{plain}      
 \newtheorem{proposition}[definition]{Proposition}
 \newtheorem{theorem}[definition]{Theorem}
 \newtheorem{corollary}[definition]{Corollary}
 \newtheorem{lemma}[definition]{Lemma}
\newtheorem*{conjecture}{Conjecture}

\newcommand{\GL}{\operatorname{GL}}

\newcommand{\R}{\mathbf{R}}
\newcommand{\Z}{\mathbf{Z}}
\newcommand{\Q}{\mathbf{Q}}

\begin{document}

\title{A survey of the Johnson homomorphisms of the automorphism groups of free groups \\ and related topics}

\author{Takao Satoh\thanks{
Work partially supported by a JSPS Research Fellowship for Young Scientists}}

\address{
Department of Mathematics, Faculty of Science Division II, \\
Tokyo University of Science\\
1-3 Kagurazaka Shinjuku-ku Tokyo 162-8601, Japan\\
e-mail:\,\tt{takao@rs.tus.ac.jp}}

\maketitle

\begin{abstract} This is a biased survey for the Johnson homomorphisms of the automorphism groups of free groups.
We just exposit some well known facts and recent developments for the Johnson homomorphisms and its related topics. 
\end{abstract}



\tableofcontents   

\section{Introduction}\label{S-int}

In the 1980s, Dennis Johnson established a new remarkable method to investigate
the group structure of the Torelli subgroup of the mapping class group of a surface
in a series of works \cite{Jo1}, \cite{Jo2}, \cite{Jo3} and \cite{Jo4}.
In particular, he gave a finite set of generators of the Torelli group, and he constructed a homomorphism $\tau$
to determine the abelianization of that group.
Today, his homomorphism $\tau$ is called the first Johnson homomorphism, and it
is generalized to Johnson homomorphisms of higher degrees.
Over the last two decades, good progress was made in the study of the Johnson homomorphisms of mapping class groups
through the works of many authors including Morita \cite{Mo1}, Hain \cite{Hai} and others.

\vspace{0.5em}

As is well known, the mapping class group of a compact oriented surface with one boundary component can be embedded into
the automorphism group of a free group by a classical work of Dehn and Nielsen in the 1910s and early 1920s.
So far, a large number of theories and research techniques to study the mapping class group have been applied to
investigate the automorphism group of a free group. The Johnson homomorphisms are one of these techniques.

\vspace{0.5em}

The definition of the Johnson homomorphisms were generalized
not only for the automorphism group of a free group but also that of any group $G$.
In this paper, although we mainly consider the case where $G$ is a free group and a free metabelian group,
we give the definition of the Johnson homomorphisms for the automorphism group of a general group $G$.
(See Section {\rmfamily \ref{S-John}}.) 

\vspace{0.5em}

We discuss the case where $G$ is a free group of finite rank,
the most basic and important case.
To put it plainly, the Johnson homomorphisms are useful tools to study the graded quotients of a certain descending filtration of
the automorphism group of a free group. To explain this, let us fix some notation.
Let $F_n$ be a free group of rank $n \geq 2$, and $\mathrm{Aut}\,F_n$ the automorphism group of $F_n$.
We denote by $H:=H_1(F_n,\Z)$ the abelianization of $F_n$. Let $\rho : \mathrm{Aut}\,F_n \rightarrow \mathrm{Aut}\,H$ be
the natural homomorphism induced from the abelianization $F_n \rightarrow H$.
The kernel of $\rho$ is called the IA-automorphism group of $F_n$, denoted by $\mathrm{IA}_n$.
The letters I and A stands for \lq\lq Identity" and \lq\lq Automorphism" respectively.
Bachmuth \cite{Ba1} called $\mathrm{IA}_n$ the IA-automorphism group since that consists of automorphisms which induce
identity automorphisms on the abelianized group $H$ of $F_n$.
The subgroup $\mathrm{IA}_n$ reflects much richness and complexity of the structure of $\mathrm{Aut}\,F_n$,
and plays important roles in various studies of $\mathrm{Aut}\,F_n$.

\vspace{0.5em}

Although the study of the IA-automorphism group has a long history
since finitely many generators were obtained for that group by Magnus \cite{Mag} in 1935,
the group structure of $\mathrm{IA}_n$ is still quite complicated.
For instance, no presentation for $\mathrm{IA}_n$ is known for $n \geq 3$.
Nielsen \cite{Ni0} showed that $\mathrm{IA}_2$ coincides with the inner automorphism group,
hence, is a free group of rank $2$.
For $n \geq 3$, however, $\mathrm{IA}_n$ is much larger than the inner automorphism group $\mathrm{Inn}\,F_n$.
Krsti\'{c} and McCool \cite{Krs} showed that $\mathrm{IA}_3$ is not finitely presentable.
For $n \geq 4$, it is not known whether $\mathrm{IA}_n$ is finitely presentable or not.

\vspace{0.5em}

Because of the complexity of the group structure of $\mathrm{IA}_n$ mentioned above,
it is sometimes not suitable to handle all of $\mathrm{IA}_n$ directly.
In order to study $\mathrm{IA}_n$ with a phased approach, we consider the Johnson filtration of $\mathrm{Aut}\,F_n$.
The Johnson filtration is one of descending central series
\[ \mathrm{IA}_n = \mathcal{A}_n(1) \supset \mathcal{A}_n(2) \supset \cdots \]
consisting of normal subgroups of $\mathrm{Aut}\,F_n$, whose first term is $\mathrm{IA}_n$.
(For details, see Sections {\rmfamily \ref{S-John}} and {\rmfamily \ref{S-Fre}}.)
The graded quotients $\mathrm{gr}^k (\mathcal{A}_n) := \mathcal{A}_n(k)/\mathcal{A}_n(k+1)$ naturally have $\mathrm{GL}(n,\Z)$-module structures,
and they are considered to be a sequence of approximations of $\mathrm{IA}_n$.

\vspace{0.5em}

To understand the graded quotients $\mathrm{gr}^k (\mathcal{A}_n)$ more closely, we use the Johnson homomorphisms
\[ \tau_k : \mathrm{gr}^k (\mathcal{A}_n) \rightarrow H^* \otimes_{\Z} \mathrm{gr}^{k+1}(\mathcal{L}_n). \]
(For the definition, see Section {\rmfamily \ref{S-John}}.)
One of the most fundamental properties of the Johnson homomorphism is that $\tau_k$ is a $\mathrm{GL}(n,\Z)$-equivariant injective homomorphism for each $k \geq 1$.
Hence, we can consider $\mathrm{gr}^k (\mathcal{A}_n)$ as a submodule of $H^* \otimes_{\Z} \mathrm{gr}^{k+1}(\mathcal{L}_n)$ whose module structure is easy to handle.
Using the Johnson homomorphisms, we can extract some valuable information about $\mathrm{IA}_n$.
For example, $\tau_1$ is just the abelianization of $\mathrm{IA}_n$ by a result due to Andreadakis \cite{And}, and $\tau_2$
was applied to determine the image of the cup product
$\cup_{\Q} : \Lambda^2 H^1(\mathrm{IA}_n,\Q) \rightarrow H^2(\mathrm{IA}_n,\Q)$ by Pettet \cite{Pet}. 

\vspace{0.5em}

In general, to determine the structure of the images, or equivalently
the cokernels, of the Johnson homomorphisms is one of the most basic problems.
In order to take advantage of representation theory, we consider the rationalization of modules.
So far, for $1 \leq k \leq 3$, the $\mathrm{GL}(n,\Q)$-module structure of the cokernel $\mathrm{Coker}(\tau_{k,\Q})$ of the rational Johnson homomorphism
$\tau_{k,\Q}:= \tau_k \otimes \mathrm{id}_{\Q}$ has been determined. (See \cite{And}, \cite{Pet} and \cite{S03} for $k=1$, $2$ and $3$ respectively,
and see also Subsection {\rmfamily \ref{Ss-JonF}}.)
In general, however, it is quite a difficult problem to determine the $\mathrm{GL}(n,\Q)$-module structure of $\mathrm{Coker}(\tau_{k,\Q})$ for arbitrary $k \geq 4$.
One reason for it is that we cannot obtain an explicit generating system of ${\mathrm{gr}}^k (\mathcal{A}_n)$ easily.

\vspace{0.5em}

To avoid this difficulty, we consider the lower central series $\Gamma_{\mathrm{IA}_n}(1)=\mathrm{IA}_n$, $\Gamma_{\mathrm{IA}_n}(2)$, $\dots$
of $\mathrm{IA}_n$.
Since the Johnson filtration is central, we have $\Gamma_{\mathrm{IA}_n}(k) \subset \mathcal{A}_n(k)$ for any $k \geq 1$.
It is conjectured that $\Gamma_{\mathrm{IA}_n}(k) = \mathcal{A}_n(k)$ for each $k \geq 1$.
The conjecture was made by Andreadakis who showed
$\Gamma_{\mathrm{IA}_2}(k) = \mathcal{A}_2(k)$ for each $k \geq 1$ and $\Gamma_{\mathrm{IA}_3}(3) = \mathcal{A}_3(3)$ in \cite{And}.
Now, it is known that $\Gamma_{\mathrm{IA}_n}(2) = \mathcal{A}_n(2)$ due to
Bachmuth \cite{Ba2} and that
$\Gamma_{\mathrm{IA}_n}(3)$ has at most finite index in $\mathcal{A}_n(3)$ due to Pettet \cite{Pet}.

\vspace{0.5em}

For each $k \geq 1$, set $\mathrm{gr}^k(\mathcal{L}_{\mathrm{IA}_n}) := \Gamma_{\mathrm{IA}_n}(k)/\Gamma_{\mathrm{IA}_n}(k+1)$.
Since $\mathrm{IA}_n$ is finitely generated as mentioned above, each $\mathrm{gr}^k(\mathcal{L}_{\mathrm{IA}_n})$ is also finitely generated
as an abelian group. Then we can define a $\mathrm{GL}(n,\Z)$-equivariant homomorphism
\[ \tau_k' : \mathrm{gr}^k(\mathcal{L}_{\mathrm{IA}_n}) \rightarrow H^* \otimes_{\Z} \mathrm{gr}^{k+1}(\mathcal{L}_n) \]
in the same way as $\tau_k$. We also call $\tau_k'$ the Johnson homomorphism of $\mathrm{Aut}\,F_n$.
Then, we can directly obtain information about the cokernel of $\tau_k'$ using finitely many generators of
$\mathrm{gr}^k(\mathcal{L}_{\mathrm{IA}_n})$.
Furthermore, if we consider the rational Johnson homomorphism $\tau_{k,\Q}' := \tau_k' \otimes \mathrm{id}_{\Q}$, then
using representation theory, we can consider $\mathrm{Coker}(\tau_{k,\Q})$ as a $\mathrm{GL}(n,\Z)$-submodule of
$\mathrm{Coker}(\tau_{k,\Q}')$. Hence, we can give an upper bound on $\mathrm{Coker}(\tau_{k,\Q})$.
In \cite{S11}, we determined the cokernel of the rational Johnson homomorphism $\tau_{k,\Q}'$ in a stable range as follows.
\begin{theorem}[Satoh, \cite{S11}] ($=$ Theorem {\rmfamily \ref{T-S11}}.)
For any $k \geq 2$ and $n \geq k+2$,
\[ \mathrm{Coker}(\tau_{k, \Q}') \cong \mathcal{C}_n^{\Q}(k). \]
\end{theorem}
Here $\mathcal{C}_n(k)$ be a quotient module of $H^{\otimes k}$ by the action of the cyclic group $\mathrm{Cyc}_k$ of order $k$ on the components: 
\[ \mathcal{C}_n(k) := H^{\otimes k} \big{/} \langle a_1 \otimes a_2 \otimes \cdots \otimes a_k - a_2 \otimes a_3 \otimes \cdots
  \otimes a_k \otimes a_1 \,|\, a_i \in H \rangle, \]
and $\mathcal{C}_n^{\Q}(k) := \mathcal{C}_n(k) \otimes_{\Z} \Q$.
The fact above is quite important if we consider an application of the study of the Johnson homomorphisms of $\mathrm{Aut}\,F_n$ to
that of the mapping class group. 

\vspace{0.5em}

Let $\mathcal{M}_{g,1}$ be the mapping class group of a compact oriented surface $\Sigma_{g,1}$
of genus $g$ with one boundary component.
By an similar way as $\tau_k$, we can define an injective $\mathrm{Sp}(2g,\Z)$-equivariant homomorphism
\[ \tau_k^{\mathcal{M}} : \mathrm{gr}^k (\mathcal{M}_{g,1}) \hookrightarrow \mathrm{Hom}_{\Z}(H, \mathrm{gr}^{k+1}(\mathcal{L}_{2g})), \]
called the $k$-th Johnson homomorphism of $\mathcal{M}_{g,1}$.
Inspired by Johnson's works, Morita has studied various aspects of the Johnson homomorphisms of the mapping class groups over the last two decades,
and gave many remarkable results. We also use $H$ for $H_1(\Sigma_{g,1}, \Z)$ by abuse of the language.
Identify $H$ with its dual module as an $\mathrm{Sp}(2g,\Z)$-module by Poincar\'{e} duality.
Let $\mathfrak{h}_{g,1}(k)$ be the kernel of a homomorphism
$H \otimes_{\Z} \mathrm{gr}^{k+1}(\mathcal{L}_{2g}) \rightarrow \mathrm{gr}^{k+2}(\mathcal{L}_{2g})$, which is induced from the Lie
bracket of the free Lie algebra.
Then Morita \cite{Mo1} showed that $\mathrm{Im}(\tau_k^{\mathcal{M}}) \subset \mathfrak{h}_{g,1}(k)$ for each $k \geq 2$.
Hence, our main problem is to determine
$\mathrm{Coker}(\tau_k^{\mathcal{M}}):=\mathfrak{h}_{g,1}(k)/\mathrm{Im}(\tau_k^{\mathcal{M}})$.
At the present stage, the Sp-module structure of the cokernel of the rational Johnson homomorphisms $\tau_{k,\Q}^{\mathcal{M}}$
are determined for $1 \leq k \leq 4$ by Johnson \cite{Jo1}, Morita \cite{Mo0}, Hain \cite{Hai} and Asada-Nakamura \cite{AN1} respectively.
Furthermore, by using trace maps, Morita \cite{Mo1} showed that the symmetric tensor product $S^k H_{\Q}$ injects into
$\mathrm{Coker}(\tau_{k,\Q}^{\mathcal{M}})$ for each odd $k \geq 3$. In 1996, in unpublished work,
Hiroaki Nakamura showed that the multiplicity of $S^k H_{\Q}$ in $\mathrm{Coker}(\tau_{k,\Q}^{\mathcal{M}})$ is one
for odd $k \geq 3$. Except for these results, there are few results for the cokernels of the Johnson homomorphisms of
the mapping class groups. Recently, based on a work of Hain \cite{Hai}, Naoya Enomoto and the author proved
\begin{theorem}[Enomoto and Satoh \cite{ES2}]($=$ Theorem {\rmfamily \ref{ES-main}}.)
For an integer $k \geq 5$ such that $k \equiv 1$ mod $4$, and $g \geq k+2$,
the irreducible Sp-module $[1^k]$ appears in $\mathrm{Coker}(\tau_{k,\Q}^{\mathcal{M}})$
with multiplicity one.
\end{theorem}
\noindent
See Section {\rmfamily \ref{S-MCG}} for details.

\vspace{0.5em}

This chapter is a survey on the Johnson homomorphisms of the automorphism group of a free group and related topics.
We just exposit some well known facts and recent developments for them without proof.
In Section {\rmfamily \ref{S-John}},
we give the definitions and some properties of the Andreadakis-Johnson filtration and the Johnson homomorphisms of the automorphism group of an
arbitrary group $G$.
In Section {\rmfamily \ref{S-Fre}}, we consider the free group case. In particular, we show that the symmetric tensor product $S^k H_{\Q}$
of $H_{\Q} := H \otimes_{\Z} \Q$ injects into
the cokernel $\mathrm{Coker}(\tau_{k,\Q})$ using the Morita trace map. Moreover, we show that $\mathrm{Coker}(\tau_{k, \Q}') \cong \mathcal{C}_n^{\Q}(k)$
for $k \geq 2$ and $n \geq k+2$.
In Section {\rmfamily \ref{S-FreM}}, we deal with the free metabelian group case. In particular, we show how the image and the cokernel of the Johnson
homomorphisms of the automorphism group of a free metabelian group were completely determined by using the Morita trace maps.
Furthermore we also show that as a Lie algebra, the abelianization of the derivation algebra of the free Lie algebra generated by $H$ was
also determined by using the Morita trace maps.

\vspace{0.5em}

In the rest of this chapter, we consider some applications of the study of the Johnson homomorphisms of the automorphism group of a free group.
In Section {\rmfamily \ref{S-BM}}, we recall the definition of the braid group and the mapping class group of a surface with boundary,
which are considered as subgroups of
the automorphism group of a free group. In this section, we give only the definition and some properties of them as a preparation for
later sections, and no topics related to the Johnson homomorphisms.
In Section {\rmfamily \ref{S-Mag}}, we show a relation between the Johnson homomorphisms and the Magnus representation obtained by Morita.
We also show the infinite generation of the kernel of the Magnus representation by using the first Johnson homomorphisms.
In Section {\rmfamily \ref{S-MCG}}, we investigate the Johnson cokernels of the mapping class group of a surface. In particular, we show that
we can detect new series $[1^k]$ in the Johnson cokernel of the mapping class group by using Theorem {\rmfamily \ref{T-S11}}.
In Section {\rmfamily \ref{S-TwC}}, we consider twisted cohomologies of the automorphism groups of a free group and a free nilpotent group.
In particular, we compute $H^1(\mathrm{Aut}\,F_n, \mathrm{IA}_n^{\mathrm{ab}})$, and show that the first Johnson homomorphism cannot be extended to
the automorphism group of a free group as a crossed homomorphism.
Moreover, for a free nilpotent group $N_{n,k}$ we show that $H^1(\mathrm{Aut}\,N_{n,k}, \Lambda^l H_{\Q})=0$ for $k \geq 3$, $n \geq k-1$ and $l \geq 3$.
Then we see that the trace map $\mathrm{Tr}_{[1^k]}$ for the exterior product $\Lambda^k H_{\Q}$ defines a non-trivial second cohomology class
in $H^2(\mathrm{Aut}\,N_{n,k}, \Lambda^k H_{\Q})$.
Finally, in Section {\rmfamily \ref{S-Cong}}, we determine the abelianization of the congruence IA-automorphism group of a free group
by using the first Johnson map defined by the Magnus expansion of a free group obtained by Kawazumi.

\section{Notation and Conventions}\label{S-Not}
\hspace*{\fill}\ 

\vspace{0.5em}

Throughout this chapter, we use the following notation and conventions.

\begin{itemize}
\item For a group $G$, the abelianization of $G$ is denoted by $G^{\mathrm{ab}}$.
\item For elements $x$ and $y$ of a group $G$, the commutator bracket $[x,y]$ of $x$ and $y$
      is defined to be $[x,y]:=xyx^{-1}y^{-1}$. For subgroups $H$ and $K$ of $G$, we denote by $[H,K]$ the commutator subgroup of $G$
      generated by $[h, k]$ for $h \in H$ and $k \in K$.
\item For a group $G$, and its quotient group $G/N$, we also denote the coset class of an element $g \in G$ by $g \in G/N$ if
      there is no confusion.

\item For a group $G$, the group $\mathrm{Aut}\,G$ acts on $G$ from the right unless otherwise noted.
      For any $\sigma \in \mathrm{Aut}\,G$ and $x \in G$,
      the action of $\sigma$ on $x$ is denoted by $x^{\sigma}$.
\item If we need to consider the left action of $\mathrm{Aut}\,G$ on $G$, we consider the usual action. Namely,
      $\sigma \cdot x := x^{\sigma^{-1}}$ for $\sigma \in \mathrm{Aut}\,G$ and $x \in G$.
      Similarly, for any right $\mathrm{Aut}\,G$-module $M$, we consider $M$ as a left $\mathrm{Aut}\,G$-module by
      $\sigma \cdot m := m^{\sigma^{-1}}$ for $\sigma \in \mathrm{Aut}\,G$ and $m \in M$.

\item For any $\Z$-module $M$, we denote $M \otimes_{\Z} \Q$ by the symbol obtained by attaching a subscript or superscript $\Q$ to $M$, like
      $M_{\Q}$ or $M^{\Q}$. Similarly, for any $\Z$-linear map $f: A \rightarrow B$,
      the induced $\Q$-linear map $f \otimes \mathrm{id}_{\Q} : A_{\Q} \rightarrow B_{\Q}$ is denoted by $f_{\Q}$
      or $f^{\Q}$.

\item For a group $G$ and a left $G$-module $M$, we set
\[\begin{split}
    \mathrm{Cros}(G,M) & := \{ f : G \rightarrow M \,|\, f : \mathrm{crossed} \,\,\, \mathrm{homomorphism} \}.
  \end{split}\]

\item For any ring $R$, we denote by $M(n, R)$ and $\mathrm{GL}(n,R)$ the ring of matrices of rank $n$ over $R$
and the general linear group over $R$ respectively.
\end{itemize}

\section{Johnson homomorphisms}\label{S-John}

In this section, we recall the definition and some properties of the Johnson homomorphisms of the automorphism group of a group.

\subsection{Graded Lie algebras associated to an N-series of a group $G$}\label{Ss-GLA}

To begin with, we recall an N-series of a group and the graded Lie algebra associated to it, based on a classical work of Lazard \cite{Laz}.

\begin{definition}
Let $G$ be a group, and $H_1, H_2, \ldots$ a sequence of subgroups of $G$ such that
\[ G = H_1 \supset H_2 \supset \cdots \supset H_k \supset \cdots \]
where $H_{k+1}$ is normal in $H_k$, and
\[ [H_k, H_l] \subset H_{k+l} \]
for any $k, l \geq 1$.
Such a sequence $\mathcal{H} := \{ H_k \}$ is called an {\it N-series} \index{N-series} of $G$. 
\end{definition}

The most familiar example of an N-series of $G$ is its lower central series $\{ \Gamma_G(k) \}$ defined by
\[ \Gamma_G(1):= G, \hspace{1em} \Gamma_G(k+1) := [\Gamma_G(k),G], \hspace{1em} k \geq 1. \]
The lower central series $\{ \Gamma_G(k) \}$ has the property that if $\{ H_k \}$ is an N-series of $G$ then $\Gamma_G(k) \subset H_k$ for any $k \geq 1$.

\vspace{0.5em}

Obviously, each of the graded quotients $\mathrm{gr}^k(\mathcal{H}) := H_k/H_{k+1}$ is an abelian group for $k \geq 1$ since $[H_k, H_k] \subset H_{2k} \subset H_{k+1}$.
We write the product in each of $\mathrm{gr}^k(\mathcal{H})$ additively.
Namely, for any $x, y \in H_k$, if we denote their coset classes modulo $H_{k+1}$ by $[x]$ and $[y]$ then
\[ [x] + [y] = [xy]. \]
Consider the direct sum
\[ \mathrm{gr}(\mathcal{H}) := \bigoplus_{k=1}^{\infty} \mathrm{gr}^k(\mathcal{H}). \]
We give the additive abelian group $\mathrm{gr}(\mathcal{H})$ a grading with the convention that the elements of $\mathrm{gr}^k(\mathcal{H})$ are homogeneous of degree $k$.

\vspace{0.5em}

For any $k$, $l \geq 1$, let us consider a bilinear alternating map
\[ [ \, , \,]_{\mathrm{Lie}} : \mathrm{gr}^k(\mathcal{H}) \times \mathrm{gr}^l(\mathcal{H}) \rightarrow \mathrm{gr}^{k+l}(\mathcal{H}) \]
defined by $[\, [x], [y] \, ]_{\mathrm{Lie}} := [\, [x, y] \,]$ for any $[x] \in \mathrm{gr}^k(\mathcal{H})$ and $[y] \in \mathrm{gr}^l(\mathcal{H})$
where $[x, y]$ is a commutator in $G$, and $[\, [x, y] \,]$ is a coset class of $[x, y]$ in $\mathrm{gr}^{k+l}(\mathcal{H})$.
Then $[ \, , \,]_{\mathrm{Lie}}$ induces a graded Lie algebra structure of the graded sum $\mathrm{gr}(\mathcal{H})$.

\begin{definition}
The Lie algebra $\mathrm{gr}(\mathcal{H})$ is called the {\it graded Lie algebra associated to the N-series} $\mathcal{H}$.
\index{graded Lie algebra associated to an N-series}
\end{definition}

This is a generalization of the well-known construction introduced by Magnus \cite{Ma2} for the lower central series of a free group.
(For example, see also Section 5.14 in \cite{MKS} or Chapter VIII in \cite{Pas} for basic materials concerning the
graded Lie algebra associated to an N-series.)

\vspace{0.5em}

For a group $G$ and its lower central series $\mathcal{L}_G := \{ \Gamma_G(k) \}$, if generators of $G$ are given, then those of the graded quotient
$\mathrm{gr}^k(\mathcal{L}_G) := \Gamma_G(k)/\Gamma_G(k+1)$ are obtained as follows.
For any $g_1, \ldots, g_k \in G$, a left-normed commutator
\[ [[ \cdots [[ g_{1},g_{i}],g_{3}], \cdots ], g_{k}] \]
of weight $k$ is called a simple $k$-fold commutator, denoted by
\[ [g_{1}, g_{2}, \cdots, g_{k} ] \]
for simplicity. Then we have

\begin{lemma}\label{l1}
If a group $G$ is generated by $g_1, \ldots , g_t$, then each of the graded quotients
$\mathrm{gr}^k(\mathcal{L}_G)$ is generated by the simple $k$-fold commutators
\[ [g_{i_1},g_{i_2}, \ldots , g_{i_k}], \hspace{1em} i_j \in \{1, \ldots , t \} \]
as an abelian group.
\end{lemma}
(For a proof see Theorem 5.4 in \cite{MKS}, for example.) This shows that if $G$ is finitely generated then so is $\mathrm{gr}^k(\mathcal{L}_G)$
for any $k \geq 1$.

\subsection{The Andreadakis-Johnson filtration of $\mathrm{Aut}\,G$}\label{Ss-AJ}

In this subsection, we consider a descending filtration of $\mathrm{Aut}\,G$.

\begin{definition}
For $k \geq 1$, the action of $\mathrm{Aut}\,G$ on each nilpotent quotient $G/\Gamma_G(k+1)$ induces a homomorphism
\[ \mathrm{Aut}\,G \rightarrow \mathrm{Aut}(G/\Gamma_G(k+1)). \]
We denote its kernel by $\mathcal{A}_G(k)$. Then the groups $\mathcal{A}_G(k)$ define a descending filtration
\[ \mathcal{A}_G(1) \supset \mathcal{A}_G(2) \supset \cdots \supset \mathcal{A}_G(k) \supset \cdots \]
of $\mathrm{Aut}\,G$. We call $\mathcal{A}_G := \{ \mathcal{A}_G(k) \}$ the {\it Andreadakis-Johnson filtration}
\index{Andreadakis-Johnson filtration} of $\mathrm{Aut}\,G$.
The first term $\mathcal{A}_G(1)$ is called the {\it IA-automorphism group}
\index{IA-automorphism group} of $G$, and is also denoted by $\mathrm{IA}(G)$.
Namely, $\mathrm{IA}(G)$ consists of automorphisms which act on the abelianization $G^{\mathrm{ab}}$ of $G$ trivially.
\end{definition}

The Andreadakis-Johnson filtration of $\mathrm{Aut}\,G$ was originally introduced by Andreadakis \cite{And} in the 1960s.
The name \lq\lq Johnson" comes from Dennis Johnson who studied the Johnson filtration and the Johnson homomorphism for the
mapping class group of a surface in the 1980s. (See Section {\rmfamily \ref{S-MCG}} for details.)
In particular, Andreadakis showed that

\begin{theorem}[Andreadakis, \cite{And}]\label{T-And} \quad
\begin{enumerate}
\item For any $k$, $l \geq 1$, $\sigma \in \mathcal{A}_G(k)$ and $x \in \Gamma_G(l)$, $x^{-1} x^{\sigma} \in \Gamma_G(k+l)$.
\item For any $k$ and $l \geq 1$, $[\mathcal{A}_G(k), \mathcal{A}_G(l)] \subset \mathcal{A}_G(k+l)$.
\item If $\displaystyle \bigcap_{k \geq 1} \Gamma_G(k) =1$, then $\displaystyle \bigcap_{k \geq 1} \mathcal{A}_G(k) =1$.
\end{enumerate}
\end{theorem}

Part (2) of Theorem {\rmfamily \ref{T-And}} shows that the Andreadakis-Johnson filtration is an N-series of $\mathrm{IA}(G)$.
Hence its $k$-th term $\mathcal{A}_G(k)$ contains $\Gamma_{\mathrm{IA}(G)}(k)$, that of the lower central series of $\mathrm{IA}(G)$.
It is a natural question to ask how different is $\mathcal{A}_G(k)$ from $\Gamma_{\mathrm{IA}(G)}(k)$. In general, however, it is quite a difficult problem
to determine whether $\mathcal{A}_G(k)$ coincides with $\Gamma_{\mathrm{IA}(G)}(k)$ or not.
The next lemma immediately follows from Theorem {\rmfamily \ref{T-And}}.

\begin{lemma}\label{L-And} \quad
\begin{enumerate}
\item For any $k$, $l \geq 1$, $\sigma \in \Gamma_{\mathrm{IA}(G)}(k)$ and $x \in \Gamma_{\mathrm{IA}(G)}(l)$, $x^{-1} x^{\sigma} \in \Gamma_G(k+l)$.
\item For any $k$ and $l \geq 1$, $[\Gamma_{\mathrm{IA}(G)}(k), \Gamma_{\mathrm{IA}(G)}(l)] \subset \Gamma_{\mathrm{IA}(G)}(k+l)$.
\item If $\displaystyle \bigcap_{k \geq 1} \Gamma_G(k) =1$, then $\displaystyle \bigcap_{k \geq 1} \Gamma_{\mathrm{IA}(G)}(k) =1$.
\end{enumerate}
\end{lemma}

Both of the graded quotients $\mathrm{gr}^k(\mathcal{A}_G) = \mathcal{A}_G(k)/\mathcal{A}_G(k+1)$ and
$\mathrm{gr}^k(\mathcal{L}_{\mathrm{IA}(G)}) = \Gamma_{\mathrm{IA}(G)}(k)/\Gamma_{\mathrm{IA}(G)}(k+1)$
are considered as sequences of approximations of $\mathrm{IA}(G)$.
To clarify their structures plays an important role in various studies of $\mathrm{IA}(G)$.

\subsection{Actions of $\mathrm{Aut}\,G/\mathrm{IA}(G)$}

In this subsection, we define actions of $\mathrm{Aut}\,G/\mathrm{IA}(G)$ on the graded quotients $\mathrm{gr}^k(\mathcal{L}_G)$, $\mathrm{gr}^k(\mathcal{A}_G)$
and $\mathrm{gr}^k(\mathcal{L}_{\mathrm{IA}(G)})$.

\vspace{0.5em}

First, since each of $\Gamma_G(k)$ is a characteristic subgroup of $G$, the group $\mathrm{Aut}\,G$ naturally acts on $\Gamma_G(k)$ from the right,
and hence on $\mathrm{gr}^k(\mathcal{L}_G) = \Gamma_G(k)/\Gamma_G(k+1)$ for any $k \geq 1$.
From Part (1) of Theorem {\rmfamily \ref{T-And}}, the restriction of this action of $\mathrm{Aut}\,G$ to $\mathrm{IA}(G)$
is trivial. Thus, an action of the quotient group $\mathrm{Aut}\,G/\mathrm{IA}(G)$ on $\mathrm{gr}^k(\mathcal{L}_G)$ is well-defined.

\vspace{0.5em}

On the other hand, since each of $\mathcal{A}_G(k)$ is a normal subgroup of $\mathrm{Aut}\,G$, the group $\mathrm{Aut}\,G$ acts on $\mathcal{A}_G(k)$ by conjugation
from the right. Namely, for any $\sigma \in \mathrm{Aut}\,G$ and $\tau \in \mathcal{A}_G(k)$, the action of $\sigma$ on $\tau$ is given by $\sigma^{-1} \tau \sigma$.
Hence, $\mathrm{Aut}\,G$ also acts on each of $\mathrm{gr}^k(\mathcal{A}_G)$ for $k \geq 1$.
From Part (2) of Theorem {\rmfamily \ref{T-And}}, the restriction of this action of $\mathrm{Aut}\,G$ to $\mathrm{IA}(G)$
is trivial. Thus, an action of the quotient group $\mathrm{Aut}\,G/\mathrm{IA}(G)$ on $\mathrm{gr}^k(\mathcal{A}_G)$ is well-defined.
Similarly, using Part (2) of Lemma {\rmfamily \ref{L-And}}, we can define an action of $\mathrm{Aut}\,G/\mathrm{IA}(G)$ on
$\mathrm{gr}^k(\mathcal{L}_{\mathrm{IA}(G)})$.

\vspace{0.5em}

Throughout the paper, we fix these three actions of $\mathrm{Aut}\,G/\mathrm{IA}(G)$ on $\mathrm{gr}^k(\mathcal{L}_G)$, $\mathrm{gr}^k(\mathcal{A}_G)$
and $\mathrm{gr}^k(\mathcal{L}_{\mathrm{IA}(G)})$.

\subsection{Johnson homomorphisms of $\mathrm{Aut}\,G$}\label{Ss-JH}

\vspace{0.5em}

To study the $\mathrm{Aut}\,G/\mathrm{IA}(G)$-module structures of the graded quotients ${\mathrm{gr}}^k (\mathcal{A}_G)$ and
$\mathrm{gr}^k(\mathcal{L}_{\mathrm{IA}(G)})$,
we use the Johnson homomorphisms of $\mathrm{Aut}\,G$.
We give their definitions and some of their properties with short proof
since it seems that there are very few articles which discuss the definition of the Johnson homomorphisms
of $\mathrm{Aut}\,G$ for any group $G$. To begin with, we prepare some lemmas.

\vspace{0.5em}

For any $\sigma \in \mathrm{Aut}\,G$ and $x \in G$, set $s_{\sigma}(x) :=x^{-1} x^{\sigma} \in G$.
By (1) of Theorem {\rmfamily \ref{T-And}, we see that if
$\sigma \in \mathcal{A}_G(k)$ and $x \in \Gamma_G(l)$, then $s_{\sigma}(x) \in \Gamma_G(k+l)$.

\begin{lemma}\label{L-ex1}
For any $\sigma$, $\tau \in \mathrm{Aut}\,G$ and $x$, $y \in G$, we have
\begin{enumerate}
\item $s_{\sigma \tau}(x) = s_{\tau}(x) \cdot s_{\sigma}(x)^{\tau} = s_{\tau}(x) s_{\sigma}(x) s_{\tau}(s_{\sigma}(x))$.
\item $s_{\sigma}(xy) = y^{-1} s_{\sigma}(x) y \cdot s_{\sigma}(y)=[y^{-1}, s_{\sigma}(x)] s_{\sigma}(x) s_{\sigma}(y)$.
\end{enumerate}
\end{lemma}
\textit{Proof.}
The equations immediately follow from
\[\begin{split}
   s_{\sigma \tau}(x) & = x^{-1} x^{\sigma \tau} = x^{-1} x^{\tau} \cdot (x^{-1} x^{\sigma})^{\tau} \\
    & = x^{-1} x^{\tau} \cdot x^{-1} x^{\sigma} \cdot (x^{-1} x^{\sigma})^{-1} \cdot (x^{-1} x^{\sigma})^{\tau}, \\
   s_{\sigma}(xy) & = y^{-1} x^{-1} x^{\sigma} y^{\sigma} = y^{-1} x^{-1} x^{\sigma} y \cdot y^{-1} y^{\sigma}.
  \end{split}\]
$\square$

\vspace{0.5em}

As a corollary to Lemma {\rmfamily \ref{L-ex1}}, we obtain

\begin{corollary}\label{C-com}
For any $\sigma$, $\tau \in \mathrm{Aut}\,G$ and $x \in G$,
\begin{enumerate}
\item $s_{\mathrm{id}}(x) = 1$, $s_{\sigma^{-1}}(x) = (s_{\sigma}(x)^{-1})^{\sigma^{-1}}$.
\item $s_{[\sigma,\tau]}(x) = (s_{\tau}(x)^{-1})^{\tau^{-1}} (s_{\sigma}(x)^{-1})^{\sigma^{-1}\tau^{-1}}$
      $s_{\tau}(x)^{\sigma^{-1}\tau^{-1}} s_{\sigma}(x)^{\tau \sigma^{-1}\tau^{-1}}$.
\end{enumerate}
\end{corollary}
\textit{Proof.}
Part (1) is trivial by the definition and (1) of Lemma {\rmfamily \ref{L-ex1}}. Part (2) is obtained from the observation
\[\begin{split}
   s_{[\sigma,\tau]}(x) & = s_{\sigma \tau (\tau \sigma)^{-1}}(x) = s_{(\tau \sigma)^{-1}}(x) \cdot s_{\sigma \tau}(x)^{\sigma^{-1} \tau^{-1}} \\
    & = (s_{\tau \sigma}(x)^{-1})^{\sigma^{-1} \tau^{-1}} (s_{\sigma \tau}(x))^{\sigma^{-1} \tau^{-1}}
  \end{split}\]
and (1) of Lemma {\rmfamily \ref{L-ex1}}. $\square$

\vspace{0.5em}

Furthermore, we see

\begin{lemma}\label{L-ex2}
For any $\sigma \in \mathcal{A}_G(k)$, $\tau \in \mathcal{A}_G(l)$ and $x \in G$,
\[ s_{[\sigma,\tau]}(x) = s_{\tau}(s_{\sigma}(x)) - s_{\sigma}(s_{\tau}(x)) \in \mathrm{gr}^{k+l+1}(\mathcal{L}_G).  \]
\end{lemma}
\textit{Proof.}
Since $s_{\sigma}(x)^{-1} \in \Gamma_G(k+1)$ and $s_{\tau}(x) \in \Gamma_G(l+1)$,
\[ s_{\sigma}(x)^{-1} s_{\tau}(x) \equiv s_{\tau}(x) s_{\sigma}(x)^{-1} \pmod{\Gamma_G(k+l+2)}, \]
and hence
\[ (s_{\sigma}(x)^{-1}s_{\tau}(x))^{\sigma^{-1}\tau^{-1}} \equiv (s_{\tau}(x)s_{\sigma}(x)^{-1})^{\sigma^{-1}\tau^{-1}} \pmod{\Gamma_G(k+l+2)}. \]
Then, from (2) of Corollary {\rmfamily \ref{C-com}}, we have
\[\begin{split}
   s_{[\sigma,\tau]}(x) & \equiv (s_{\tau}(x)^{-1})^{\tau^{-1}} s_{\tau}(x)^{\sigma^{-1} \tau^{-1}} (s_{\sigma}(x)^{-1})^{\sigma^{-1} \tau^{-1}}
     s_{\sigma}(x)^{\tau \sigma^{-1} \tau^{-1}} \\
   & = (s_{\sigma^{-1}}(s_{\tau}(x)))^{\tau^{-1}} \cdot (s_{\tau}(s_{\sigma}(x)))^{\sigma^{-1} \tau^{-1}} \\
   & = \big{(} (s_{\sigma}(s_{\tau}(x)))^{-1} \big{)}^{\sigma^{-1} \tau^{-1}} \cdot (s_{\tau}(s_{\sigma}(x)))^{\sigma^{-1} \tau^{-1}}
  \end{split}\]
modulo $\Gamma_G(k+l+2)$. Since $\sigma$ and $\tau \in \mathrm{IA}(G)$ act on $\mathrm{gr}^{k+l+1}(\mathcal{L}_G)$ trivially,
we obtain the required result. $\square$

\vspace{0.5em}

Now, for any $\sigma \in \mathcal{A}_G(k)$, consider a map
$\tilde{\tau}_k(\sigma) :G^{\mathrm{ab}} \rightarrow \mathrm{gr}^{k+1}({\mathcal{L}}_{G})$ defined by
\[ x \pmod{\Gamma_G(2)} \hspace{1em} \mapsto \hspace{1em} s_{\sigma}(x) \pmod{\Gamma_G(k+2)} \]
for any $x \in G$. This is well-defined. In fact, if $x \equiv y \pmod{\Gamma_G(2)}$, there exists some $c \in \Gamma_G(2)$ such that
$y=xc$. Then, by (2) of Lemma {\rmfamily \ref{L-ex1}}, we see
\[ s_{\sigma}(y) = s_{\sigma}(xc) =[c^{-1}, s_{\sigma}(x)] s_{\sigma}(x) s_{\sigma}(c) \equiv s_{\sigma}(x) \pmod{\Gamma_G(k+2)}. \]
Similarly, using (2) of Lemma {\rmfamily \ref{L-ex1}}, we see that $\tilde{\tau}_k(\sigma)$ is a homomorphism between abelian groups.

\vspace{0.5em}

Furthermore, a map $\tilde{\tau}_k : \mathcal{A}_G(k) \rightarrow \mathrm{Hom}_{\Z}(G^{\mathrm{ab}}, \mathrm{gr}^{k+1}(\mathcal{L}_G))$ defined by
\[ \sigma \hspace{0.5em} \mapsto \hspace{0.5em} \tilde{\tau}_k(\sigma) \]
is a homomorphism. In fact, for any $\sigma, \tau \in \mathcal{A}_G(k)$, by (1) of Lemma {\rmfamily \ref{L-ex1}} we have
\[\begin{split}
   \tilde{\tau}_k(\sigma \tau)(x) & = s_{\sigma \tau}(x) = s_{\tau}(x) s_{\sigma}(x) s_{\tau}(s_{\sigma}(x)) \\
     & = s_{\tau}(x) + s_{\sigma}(x) = \tilde{\tau}_k(\sigma)(x) + \tilde{\tau}_k(\tau)(x)
  \end{split}\]
for any $x \in G^{\mathrm{ab}}$. From the definition, it is easy to see that the kernel of $\tilde{\tau}_k$ is just $\mathcal{A}_G(k+1)$.
Therefore $\tilde{\tau}_k$ induces an injective homomorphism
\[ \tau_k : \mathrm{gr}^k (\mathcal{A}_G) \hookrightarrow \mathrm{Hom}_{\Z}(G^{\mathrm{ab}}, \mathrm{gr}^{k+1}(\mathcal{L}_G)). \]

\begin{definition}
For each $k \geq 1$, we call the homomorphisms $\tilde{\tau}_k$ and $\tau_k$ the $k$-th {\it Johnson homomorphisms}
\index{Johnson homomorphisms} of $\mathrm{Aut}\,G$.
\end{definition}

Here we remark

\begin{lemma}
For each $k \geq 1$, $\tau_k$ is an $\mathrm{Aut}\,G/\mathrm{IA}(G)$-equivariant homomorphism.
\end{lemma}
\textit{Proof.}
For any $\sigma \in \mathrm{Aut}\,G$ and $\tau \in \mathcal{A}_G(k)$, we see
\[\begin{split}
   \tau_k(\tau \cdot \sigma)(x) & = \tau_k(\sigma^{-1} \tau \sigma)(x) = s_{\sigma^{-1} \tau \sigma}(x), \\
   (\tau_k(\tau) \cdot \sigma)(x) & = ( \tau_k(\tau)(x^{\sigma^{-1}}) )^{\sigma} = s_{\tau}(x^{\sigma^{-1}})^{\sigma} = ((x^{\sigma^{-1}})^{-1} x^{\sigma^{-1}\tau})^{\sigma} \\
               & = s_{\sigma^{-1} \tau \sigma}(x)
  \end{split}\]
for any $x \in G^{\mathrm{ab}}$. Hence we have $\tau_k(\tau \cdot \sigma) = \tau_k(\tau) \cdot \sigma$. This means $\tau_k$ is an $\mathrm{Aut}\,G$-equivariant
homomorphism. Since $\mathrm{IA}(G)$ acts on $\mathrm{gr}^k (\mathcal{A}_G)$ and $\mathrm{Hom}_{\Z}(G^{\mathrm{ab}}, \mathcal{L}_G(k+1))$ trivially,
we obtain the required result. $\square$

\vspace{0.5em}

From the viewpoint of the study of the $\mathrm{Aut}\,G/\mathrm{IA}(G)$-module structure of $\mathrm{gr}^k(\mathcal{A}_G)$,
it is a natural and basic problem to determine the image of $\tau_k$. In general, this is difficult even in the case where $G$ is a free group.
(See Section {\rmfamily \ref{S-Fre}}.)

\vspace{1em}

Now, we recall the derivation algebra of a Lie algebra. (See Section 8 of Chapter II in \cite{Bou} for details on the derivation algebra.)

\begin{definition}
Let $\mathfrak{g}$ be a Lie algebra over $\Z$.
A $\Z$-linear map $f : \mathfrak{g} \rightarrow \mathfrak{g}$ is called a derivation of $\mathfrak{g}$ if
$f$ satisfies
\[ f([a,b]) = [f(a),b]+ [a,f(b)] \]
for any $a, b \in \mathfrak{g}$.
Let $\mathrm{Der}(\mathfrak{g})$ be a set of all derivations of $\mathfrak{g}$. Then $\mathrm{Der}(\mathfrak{g})$ has a Lie algebra structure
over $\Z$. This Lie algebra $\mathrm{Der}(\mathfrak{g})$ is called the {\it derivation algebra}
\index{derivation algebra} of $\mathfrak{g}$.
\end{definition}

\vspace{0.5em}

We define a grading of $\mathrm{Der}(\mathrm{gr}(\mathcal{L}_G))$ as follows.
For $k \geq 0$, the degree $k$ part of $\mathrm{Der}(\mathcal{L}_G)$ is defined to be
\[ \mathrm{Der}(\mathrm{gr}(\mathcal{L}_G))(k) := \{ f \in \mathrm{Der}(\mathrm{gr}(\mathcal{L}_G)) \,|\, f(a)
     \in \mathrm{gr}^{k+1}(\mathcal{L}_G), \,\,\, a \in G^{\mathrm{ab}} \}. \]
Then, we have
\[ \mathrm{Der}(\mathrm{gr}(\mathcal{L}_G)) = \bigoplus_{k \geq 0} \mathrm{Der}(\mathrm{gr}(\mathcal{L}_G))(k) \]
as an $\mathrm{Aut}\,G/\mathrm{IA}(G)$-module. For any derivation $f$ of $\mathrm{gr}(\mathcal{L}_G)$, since the image of $f$ is completely determined by
that of the degree $1$ part $\mathrm{gr}^1(\mathcal{L}_G)=G^{\mathrm{ab}}$, the restriction of an element of $\mathrm{Der}(\mathrm{gr}(\mathcal{L}_G))(k)$
to $\mathrm{gr}^1(\mathcal{L}_G)$ induces an $\mathrm{Aut}\,G/\mathrm{IA}(G)$-equivariant embedding from $\mathrm{Der}(\mathrm{gr}(\mathcal{L}_G))(k)$
into $\mathrm{Hom}_{\Z}(G^{\mathrm{ab}}, \mathrm{gr}^{k+1}(\mathcal{L}_G))$. In this paper, by this embedding,
we consider $\mathrm{Der}(\mathrm{gr}(\mathcal{L}_G))(k)$
as an $\mathrm{Aut}\,G/\mathrm{IA}(G)$-submodule of $\mathrm{Hom}_{\Z}(G^{\mathrm{ab}}, \mathrm{gr}^{k+1}(\mathcal{L}_G))$
for each $k \geq 0$.

\vspace{0.5em}

Let $\mathrm{Der}^+(\mathrm{gr}(\mathcal{L}_G))$ be a graded Lie subalgebra of $\mathrm{Der}(\mathrm{gr}(\mathcal{L}_G))$ with positive degree.
Consider a direct sum of the Johnson homomorphisms
\[ \bigoplus_{k \geq 1} \tau_k : \mathrm{gr}(\mathcal{A}_G) \hookrightarrow \mathrm{Der}^+(\mathrm{gr}(\mathcal{L}_G)). \]
By Lemma {\rmfamily \ref{L-ex2}}, the map $\bigoplus_{k \geq 1} \tau_k$ is a Lie algebra homomorphism. This is called the
{\it total Johnson homomorphism} \index{total Johnson homomorphism} \index{Johnson homomorphism!total} of $\mathrm{Aut}\,G$.

\vspace{1em}

To end this section, we define the Johnson homomorphisms for the graded quotients $\mathrm{gr}^k(\mathcal{L}_{\mathrm{IA}(G)})$.

\begin{definition}
The restriction of the homomorphism $\tilde{\tau}_k : \mathcal{A}_G(k) \rightarrow \mathrm{Hom}_{\Z}(G^{\mathrm{ab}}$, $\mathrm{gr}^{k+1}(\mathcal{L}_G))$
to $\Gamma_{\mathrm{IA}(G)}(k)$ induces a homomorphism
\[ \tau_k' : \mathrm{gr}^k (\mathcal{L}_{\mathrm{IA}(G)}) \rightarrow \mathrm{Hom}_{\Z}(G^{\mathrm{ab}}, \mathrm{gr}^{k+1}(\mathcal{L}_G)). \]
We also call $\tau_k'$ the $k$-th {\it Johnson homomorphism} of $\mathrm{Aut}\,G$.
\end{definition}

We remark that each of $\tau_k'$ is an $\mathrm{Aut}\,G/\mathrm{IA}(G)$-equivariant homomorphism, by the same argument as for $\tau_k$.
The Andreadakis-Johnson filtration $\mathcal{A}_G$ coincides with the lower central series $\mathcal{L}_{\mathrm{IA}(G)}$
of $\mathrm{IA}(G)$ if and only if all of the Johnson homomorphisms $\tau_k'$ are injective.

\vspace{0.5em}

By an argument similar to that of $\bigoplus_{k \geq 1} \tau_k$, we see that a direct sum
\[ \bigoplus_{k \geq 1} \tau_k' : \mathrm{gr}(\mathcal{L}_{\mathrm{IA}(G)}) \rightarrow \mathrm{Der}^+(\mathrm{gr}(\mathcal{L}_G)) \]
of $\tau_k'$ is an $\mathrm{Aut}\,G/\mathrm{IA}(G)$-equivariant Lie algebra homomorphism.
We also call it the {\it total Johnson homomorphism} of $\mathrm{Aut}\,G$.

\section{Free group case}\label{S-Fre}

In this section we apply the argument in Section {\rmfamily \ref{S-John}}
to the case where $G$ is a free group $F_n$ of rank $n$ with a basis $x_1, \ldots, x_n$.

\vspace{0.5em}

First, we denote the abelianization of $F_n$ by $H$, and its dual group by $H^* :=\mathrm{Hom}_{\Z}(H,\Z)$.
If we fix the basis of $H$ as a free abelian group induced from the basis $x_1, \ldots , x_n$ of $F_n$, we can identify
$\mathrm{Aut}\,F_n^{\mathrm{ab}}=\mathrm{Aut}(H)$ with the general linear group $\mathrm{GL}(n,\Z)$.
By a classical work due to Nielsen, it is well-known that the homomorphism $\rho : \mathrm{Aut}\,F_n \rightarrow \mathrm{GL}(n,\Z)$ induced from
the abelianization of $F_n$ is surjective. (See Proposition 4.4 in \cite{Lyn}, for example.)
Hence we can identify $\mathrm{Aut}(H)/\mathrm{IA}(F_n)$ with $\mathrm{GL}(n,\Z)$.
Namely, we have a group extension
\[ 1 \rightarrow \mathrm{IA}(F_n) \rightarrow \mathrm{Aut}\,F_n \rightarrow \mathrm{GL}(n,\Z) \rightarrow 1. \]

\vspace{0.5em}

In the following, for simplicity, we omit the capital \lq\lq $F$" from a symbol attached with $F_n$ if there is no confusion.
For example, we write $\Gamma_n(k)$, $\mathrm{gr}^k(\mathcal{L}_n)$, $\mathrm{IA}_n$, $\mathcal{A}_n(k)$ and $\mathrm{gr}^k(\mathcal{A}_n)$
for $\Gamma_{F_n}(k)$, $\mathrm{gr}^k(\mathcal{L}_{F_n})$, $\mathrm{IA}(F_n)$ and $\mathrm{gr}^k(\mathcal{A}_{F_n})$ respectively.

\vspace{0.5em}

Although our ultimate purpose is to clarify the $\mathrm{GL}(n,\Z)$-module structures of $\mathrm{gr}^k(\mathcal{A}_n)$,
or equivalently those of $\mathrm{Coker}(\tau_k)$,
to take advantage of the representation theory for the general linear group $\mathrm{GL}(n,\Q)$, we consider the rationalization as follows.
In general, $\mathrm{GL}(n,\Z)$ can be considered as a subgroup of $\mathrm{GL}(n,\Q)$ in a natural way. Then the actions of $\mathrm{GL}(n,\Z)$
on $\mathrm{gr}_{\Q}^k(\mathcal{A}_n)$, $H_{\Q}^* \otimes_{\Q} \mathrm{gr}_{\Q}^{k+1}(\mathcal{L}_n)$ and $\mathrm{Coker}(\tau_{k,\Q})$ naturally
extend to those of $\mathrm{GL}(n,\Q)$. 

\vspace{0.5em}

We sum up several basic problems for the study of the Johnson homomorphisms of $\mathrm{Aut}\,F_n$.
\begin{enumerate}
\item Determine the $\mathrm{GL}(n,\Q)$-irreducible decompositions of
\[ \mathrm{Im}(\tau_{k,\Q}) \cong \mathrm{gr}_{\Q}^k(\mathcal{A}_n) \hspace{1em} \mathrm{and} \hspace{1em} \mathrm{Coker}(\tau_{k,\Q}). \]
      Furthermore, give a formula for the dimension of $\mathrm{gr}_{\Q}^k(\mathcal{A}_n)$, which is equal to
      the rank of $\mathrm{gr}^k(\mathcal{A}_n)$.
\item ({\it Andreadakis's conjecture}) \index{Andreadakis's conjecture}
      Determine whether $\mathcal{A}_n(k)$ coincides with $\Gamma_{\mathrm{IA}_n}(k)$ or not, for $n \geq 3$ and
      $k \geq 2$.
\item Determine whether each of $\mathcal{A}_n(k)$ is finitely generated or not, for $k \geq 2$. Moreover, what about $H_1(\mathcal{A}_n(k),\Z)$?
\end{enumerate}

\subsection{Free Lie algebras and their derivations}\label{Ss-FLA}

By a classical work of Magnus, it is known that each of $\mathrm{gr}^k(\mathcal{L}_n)$ is a free abelian group, and that
the graded Lie algebra $\mathrm{gr}(\mathcal{L}_n)$ is isomorphic to the free Lie algebra generated by $H$.
(For example, see \cite{MKS} or \cite{Reu} for details for a free Lie algebra.)
Witt \cite{Wit} gave the rank of $\mathrm{gr}^k(\mathcal{L}_n)$ by
\begin{equation}\label{ex-witt}
 \mathrm{rank}_{\Z}(\mathrm{gr}^k(\mathcal{L}_n))=\frac{1}{k} \sum_{d | k} \mu(d) n^{\frac{k}{d}}
\end{equation}
where $\mu$ is the M$\ddot{\mathrm{o}}$bius function.

\vspace{0.5em}

The graded Lie algebra $\mathrm{gr}(\mathcal{L}_n)$ is considered as a Lie subalgebra of the tensor algebra generated by $H$ as follows.
Let
\[ T(H):= \Z \oplus H \oplus H^{\otimes 2} \oplus \cdots \]
be the tensor algebra of $H$ over $\Z$. Then $T(H)$ is the
universal enveloping algebra of the free Lie algebra $\mathrm{gr}(\mathcal{L}_n)$, and the natural map
$\iota : \mathrm{gr}(\mathcal{L}_n) \rightarrow T(H)$ defined by
\[ [X,Y] \mapsto X \otimes Y - Y \otimes X \]
for $X$, $Y \in \mathrm{gr}(\mathcal{L}_n)$ 
is an injective graded Lie algebra homomorphism by the Poincar\'{e}-Birkoff-Witt's theorem.
We denote by $\iota_k$ the homomorphism which is
the degree $k$ part of $\iota$, and
consider $\mathrm{gr}^k(\mathcal{L}_n)$ as a $\mathrm{GL}(n,\Z)$-submodule $H^{\otimes k}$ through $\iota_k$.

\vspace{0.5em}

Consider the derivation algebra $\mathrm{Der}(\mathrm{gr}(\mathcal{L}_n))$ of the free Lie algebra $\mathrm{gr}(\mathcal{L}_n)$.
By the universality of the free Lie algebra $\mathrm{gr}(\mathcal{L}_n)$, the embedding
\[ \mathrm{Der}(\mathrm{gr}(\mathcal{L}_n))(k) \hookrightarrow H^* {\otimes}_{\Z} \mathrm{gr}^{k+1}(\mathcal{L}_n) \]
as mentioned above is surjective.
Namely, the degree $k$-part $\mathrm{Der}(\mathrm{gr}(\mathcal{L}_n))(k)$ is considered as
\[ \mathrm{Hom}_{\Z}(H,\mathrm{gr}^{k+1}(\mathcal{L}_n)) = H^* {\otimes}_{\Z} \mathrm{gr}^{k+1}(\mathcal{L}_n) \]
for each $k \geq 1$. (See Section 8 of Chapter II in \cite{Bou} for details.)

\subsection{The IA-automorphism group of $F_n$}\label{Ss-IAF}

Now, we consider the IA-automorphism group of $F_n$.
It is known by work of Nielsen \cite{Ni0} that $\mathrm{IA}_2$ coincides with the inner automorphsim group $\mathrm{Inn}\,F_2$
of $F_2$. Namely, $\mathrm{IA}_2$ is a free group of rank $2$.
However, $\mathrm{IA}_n$ for $n \geq 3$ is much larger than the inner automorphism group $\mathrm{Inn}\,F_n$.
In fact, Magnus \cite{Mag} showed
\begin{theorem}[Magnus, \cite{Mag}]\label{T-Mag}
 for any $n \geq 3$, the IA-automorphism group $\mathrm{IA}_n$ of $F_n$ is finitely generated by automorphisms
\[ K_{ij} : \begin{cases}
               x_i &\mapsto {x_j}^{-1} x_i x_j, \\
               x_t &\mapsto x_t, \hspace{4em} (t \neq i)
              \end{cases}\]
for distinct $i$, $j \in \{ 1, 2, \ldots , n \}$ and
\[  K_{ijk} : \begin{cases}
               x_i &\mapsto x_i x_j x_k {x_j}^{-1} {x_k}^{-1},  \\
               x_t &\mapsto x_t, \hspace{4em} (t \neq i)
              \end{cases}\] 
for distinct $i$, $j$, $k \in \{ 1, 2, \ldots , n \}$ such that $j<k$.
\end{theorem}

\vspace{0.5em}

For any $n \geq 3$, although a generating set of $\mathrm{IA}_n$ is well known as above, no presentation for
$\mathrm{IA}_n$ is known.
For $n=3$, Krsti\'{c} and McCool \cite{Krs} showed that $\mathrm{IA}_3$ is not finitely presentable.
For $n \geq 4$, it is also not known whether $\mathrm{IA}_n$ is finitely presentable or not.

\subsection{Johnson homomorphisms of $\mathrm{Aut}\,F_n$}\label{Ss-JonF}

In this subsection, we sum up some results for the Johnson homomorphisms of $\mathrm{Aut}\,F_n$.

\vspace{0.5em}

In the following, we denote by $S^k M$ and $\Lambda^k M$
the symmetric tensor product and the exterior product of a $\Z$-module $M$ of degree $k$ respectively.

\subsubsection{$\tau_1$} \quad

\vspace{0.5em}

Andreadakis \cite{And} calculated the images of the Magnus generators of $\mathrm{IA}_n$ by the first Johnson homomorphism as
\begin{equation}\label{imM}
 \tau_1(K_{ij})= x_i^* \otimes [x_i,x_j], \hspace{1em} \tau_1(K_{ijk})= x_i^* \otimes [x_j,x_k],
\end{equation}
and showed that $\tau_1$ is surjective.

\vspace{0.5em}

Furthermore, recently, Cohen-Pakianathan \cite{Co1, Co2}CFarb \cite{Far} and Kawazumi \cite{Kaw} independently showed that
$\tau_1$ induces the abelianization of $\mathrm{IA}_n$. Namely, for any $n \geq 3$, we have
\begin{equation}\label{CPFK}
\mathrm{IA}_n^{\mathrm{ab}} \cong H^* \otimes_{\Z} \Lambda^2 H
\end{equation}
as a $\mathrm{GL}(n,\Z)$-module.

\vspace{1em}

Here we remark that using (\ref{CPFK}), we can recursively calculate $\tau_k(\sigma)=\tau_k'(\sigma)$ for any $\sigma \in \Gamma_{\mathrm{IA}_n}(k)$.
These computations are perhaps most easily explained with some examples, so we give two here.
For the standard basis $x_1, \ldots , x_n$ of $H$ induced from
the generators of $F_n$, let $x_1^*, \ldots , x_n^*$ be the dual basis of $H^*$.
For distinct $1 \leq i, j, k, l \leq n$, we have
\[\begin{split}
  \tau_2([K_{ij},K_{jik}]) &= x_i^* \otimes ([s_{K_{jik}}(x_i), x_j] + [x_i, s_{K_{jik}}(x_j)]) \\
                            & \hspace{4em} - x_j^* \otimes ([s_{K_{ij}}(x_i),x_k] + [x_i, s_{K_{ij}}(x_k)] ), \\
                            &= x_i^* \otimes [x_i,[x_i,x_k]] - x_j^* \otimes [[x_i,x_j],x_k]
  \end{split}\] 
and
\[\begin{split}
  \tau_3 & ([K_{ij}, K_{jik},K_{il}]) \\
                                   &= x_i^* \otimes ([s_{K_{il}}(x_i),[x_i,x_k]] + [x_i,[s_{K_{il}}(x_i),x_k]]
                                        + [x_i,[x_i,s_{K_{il}}(x_k)]]), \\
                                   & \hspace{2em} - x_j^* \otimes ([[s_{K_{il}}(x_i),x_j],x_k] + [[x_i,s_{K_{il}}(x_j)],x_k]
                                        + [[x_i,x_j],s_{K_{il}}(x_k)]) \\
                                   & \hspace{2em} - x_i^* \otimes ([s_{[K_{ij},K_{jik}]}(x_i),x_l]
                                        + [x_i, s_{[K_{ij},K_{jik}]}(x_l)] ) \\
                                   &= x_i^* \otimes [[x_i,x_l],[x_i,x_k]] + x_i^* \otimes [x_i,[[x_i,x_l],x_k]] - x_j^* \otimes [[[x_i,x_l],x_j],x_k] \\
                                   &  \hspace{2em} - x_i^* \otimes [[x_i,[x_i,x_k]],x_l].
  \end{split}\]

\subsubsection{$\tau_2$ and $\tau_{3,\Q}$} \quad

\vspace{0.5em}

To begin with, we mention the difference between the Johnson filtration $\{ \mathcal{A}_n(k) \}$ and the lower central series $\{ \Gamma_{\mathrm{IA}_n}(k) \}$
of $\mathrm{IA}_n$.

\vspace{0.5em}

In \cite{And}, Andreadakis showed that $\mathcal{A}_2(k) = \Gamma_{\mathrm{IA}_2}(k)$ for any $k \geq 1$, and that $\mathcal{A}_3(3) = \Gamma_{\mathrm{IA}_3}(3)$.
Since $\tau_1$ is the abelianization of $\mathrm{IA}_n$ as mentioned above, we see that $\mathcal{A}_n(2) = \Gamma_{\mathrm{IA}_n}(2)$ for any $n \geq 2$.
This fact was originally shown by Bachmuth in \cite{Ba2}.
On the other hand, Pettet \cite{Pet} showed that $\Gamma_{\mathrm{IA}_n}(3)$ has at most finite index in $\mathcal{A}_n(3)$ for $n \geq 3$.
\begin{conjecture}[Andreadakis]\label{Conj-And}
$\mathcal{A}_n(k) = \Gamma_{\mathrm{IA}_n}(k)$ for any $k \geq 3$.
\end{conjecture}
There are, however, few results for the difference between $\mathcal{A}_n(k)$ and $\Gamma_{\mathrm{IA}_n}(k)$.

\vspace{0.5em}

Based on the above facts, we can directly calculate the image of $\tau_2$ and $\tau_{3,\Q}$.
Pettet \cite{Pet} showed
\begin{theorem}[Pettet, \cite{Pet}]\label{T-Pet} For any $n \geq 3$,
\[ \mathrm{Coker}(\tau_{2,\Q}) \cong S^2 H_{\Q}. \]
\end{theorem}
More precisely, she gave the $\mathrm{GL}(n,\Q)$-decomposition of $\mathrm{gr}_{\Q}^2(\mathcal{A}_n)$.
On the other hand, it is also known by work of Morita, mentioned in Subsection
{\rmfamily \ref{Ss-Mor}} below, that $S^2 H_{\Q}$ appears in $\mathrm{Coker}(\tau_{2,\Q})$. Then combining these two facts, we obtain Theorem
{\rmfamily \ref{T-Pet}}. In \cite{S03}, we show that the isomorphism above holds over $\Z$. Namely, we have
$\mathrm{Coker}(\tau_{2}) \cong S^2 H$.

\vspace{0.5em}

In addition to the above, Pettet \cite{Pet} also gave the $\mathrm{GL}(n,\Q)$-decomposition of the image of the rational cup product
\[ \cup : \Lambda^2 H^1(\mathrm{IA}_n,\Q) \rightarrow H^2(\mathrm{IA}_n, \Q). \]
We remark that it is not determined whether $H^2(\mathrm{IA}_n,\Q)$ coincides with $\mathrm{Im}(\cup)$ or not.

\vspace{0.5em}

In \cite{S03}, using the fact that $\Gamma_{\mathrm{IA}_n}(3)$ has at most finite index in $\mathcal{A}_n(3)$, by work of Pettet
mentioned above, we showed
\begin{theorem}[Satoh, \cite{S03}] For any $n \geq 3$,
\[ \mathrm{Coker}(\tau_{3,\Q}) \cong S^3 H_{\Q} \oplus \Lambda^3 H_{\Q}. \]
\end{theorem}
In general, for any $n \geq 3$ and $k \geq 4$, the $\mathrm{GL}(n,\Q)$-module structures of ${\mathrm{gr}}_{\Q}^k(\mathcal{A}_n)$ and
$\mathrm{Coker}(\tau_{k,\Q})$ are not yet determined.

\subsection{Morita's trace maps}\label{Ss-Mor} \quad

\vspace{0.5em}

In the calculation of the cokernel of $\tau_{k,\Q}$ for $k=2$ and $3$ above, we see that $S^k H_{\Q}$ appears in its irreducible decomposition.
It is known that $S^k H_{\Q}$ appears in the irreducible decomposition of $\mathrm{Coker}(\tau_{k,\Q})$ for any $k \geq 2$, due to Morita.
In this subsection, we review Morita's pioneer work to study $\mathrm{Coker}(\tau_{k,\Q})$ for a general $k \geq 2$.

\vspace{0.5em}

For each $k \geq 1$,
let $\varphi^k : H^* {\otimes}_{\Z} H^{\otimes (k+1)} \rightarrow H^{\otimes k}$
be the contraction map with respect to the first component, defined by
\[ x_i^* \otimes x_{j_1} \otimes \cdots \otimes x_{j_{k+1}} \mapsto x_i^*(x_{j_1}) \, \cdot
    x_{j_2} \otimes \cdots \otimes \cdots \otimes x_{j_{k+1}}. \]
Using the natural embedding ${\iota}_{k+1} : \mathrm{gr}^{k+1}(\mathcal{L}_n) \rightarrow H^{\otimes (k+1)}$,
we construct a $\mathrm{GL}(n,\Z)$-equivariant homomorphism
\[ \Phi^k = {\varphi}^{k} \circ ({id}_{H^*} \otimes {\iota}_{k+1})
    : H^* {\otimes}_{\Z} \mathrm{gr}^{k+1}(\mathcal{L}_n) \rightarrow H^{\otimes k}. \]
We also call $\Phi^k$ a contraction map.

\begin{definition}
For each $k \geq 2$, the composition map of $\Phi^k$ and the natural projection $H^{\otimes k} \rightarrow S^k H$ is called
{\it Morita's trace map} \index{Morita's trace map} of degree $k$,
denoted by
\[ \mathrm{Tr}_{[k]} : H^* {\otimes}_{\Z} \mathrm{gr}^{k+1}(\mathcal{L}_n) \rightarrow S^k H. \]
\end{definition}

We remark that Morita \cite{Mo1} originally defined Morita's trace using the Magnus representation of $\mathrm{Aut}\,F_n$, and that
Morita's original definition is equivalent to the above. (See Subsection {\rmfamily \ref{Ss-MJ}}.)
Morita proved

\begin{theorem}[Morita, \cite{Mo3}]
For any $n \geq 3$ and $k \geq 2$,
\begin{enumerate}
\item $\mathrm{Tr}_{[k]}$ is surjective.
\item $\mathrm{Tr}_{[k]} \circ \tau_k \equiv 0$.
\end{enumerate}
\end{theorem}

As a corollary, we obtain

\begin{corollary}\label{C-Mor}
For any $n \geq 3$ and $k \geq 2$, $S^k H_{\Q}$ appears in the irreducible decomposition of $\mathrm{Coker}(\tau_{k,\Q})$.
\end{corollary}

We call each of $S^k H_{\Q}$ the {\it Morita obstruction}. \index{Morita obstruction}
Here the term \lq\lq obstruction" means an obstruction for the surjectivity of the Johnson
homomorphism $\tau_k$.
At the present stage, there are few results on irreducible components of $\mathrm{Coker}(\tau_{k,\Q})$
for general $k \geq 2$ except for the Morita obstruction.

\vspace{0.5em}

We remark that recently we \cite{ES1} showed that the multiplicity of the Morita obstruction $S^k H_{\Q}$ in $\mathrm{Coker}(\tau_{k,\Q})$ is one.
(See also Theorem {\rmfamily \ref{T-ES1}}.)

\subsection{Cokernels of $\tau_k'$}\label{Ss-Low} \quad

\vspace{0.5em}

The goal of this subsection is to give an upper bound on $\mathrm{Coker}(\tau_{k,\Q})$ as a $\mathrm{GL}(n,\Q)$-module.
In order to do this, we consider the Johnson homomorphism
\[ \tau_k' : \mathrm{gr}^k(\mathcal{L}_{\mathrm{IA}_n}) \rightarrow H^* \otimes_{\Z} \mathrm{gr}^{k+1}(\mathcal{L}_n). \]
Since $\tau_k'$ factors through $\tau_k$, we see that $\mathrm{Im}(\tau_k') \subset \mathrm{Im}(\tau_k)$.
Hence using the representation theory for $\mathrm{GL}(n,\Q)$, we can consider $\mathrm{Coker}(\tau_{k,\Q})$ as a $\mathrm{GL}(n,\Q)$-submodule of
$\mathrm{Coker}(\tau_{k,\Q}')$. This means that $\mathrm{Coker}(\tau_{k,\Q}')$ can be regarded as an upper bound on $\mathrm{Coker}(\tau_{k,\Q})$.

\vspace{0.5em}

There are at least three reasons to study $\mathrm{Coker}(\tau_{k,\Q}')$. The first one is that
we can directly compute $\mathrm{Coker}(\tau_{k,\Q}')$ using finitely many generators of $\mathrm{gr}_{\Q}^k(\mathcal{L}_{\mathrm{IA}_n})$
induced from the Magnus generators of $\mathrm{IA}_n$ and Lemma {\rmfamily \ref{l1}}.
The second one is that this research would be applied to the study the Andreadakis conjecture.
If the conjecture is true, we get $\mathrm{Coker}(\tau_{k}')=\mathrm{Coker}(\tau_{k})$ for all $k \geq 1$.
Hence it seems to be important to determine the $\mathrm{GL}(n,\Q)$-module structure of $\mathrm{Coker}(\tau_{k, \Q}')$
on the study of the difference between $\{ \mathcal{A}_n(k) \}$ and $\{ \Gamma_{\mathrm{IA}_n}(k) \}$.
The final one is the most interesting motivation for a topological application.
By using the structure of $\mathrm{Coker}(\tau_{k,\Q}')$, we can study the cokernel of the
Johnson homomorphisms of the mapping class group of a surface. This is discussed more precisely in Section {\rmfamily \ref{S-MCG}} later.

\vspace{1em}

For $1 \leq k \leq 3$, we have $\mathrm{Coker}(\tau_{k,\Q}')=\mathrm{Coker}(\tau_{k,\Q})$, and hence they have been completely determined.
In \cite{S09}, we give the irreducible decomposition of $\mathrm{Coker}(\tau_{4,\Q}')$ for $n \geq 6$.
Furthermore, in \cite{S11}, we determined $\mathrm{Coker}(\tau_{k,\Q}')$ in a stable range as follows.
Let $\mathcal{C}_n(k)$ be a quotient module of $H^{\otimes k}$ by the action of the cyclic group $\mathrm{Cyc}_k$ of order $k$ on the components: 
\[ \mathcal{C}_n(k) := H^{\otimes k} \big{/} \langle a_1 \otimes a_2 \otimes \cdots \otimes a_k - a_2 \otimes a_3 \otimes \cdots \otimes a_k \otimes a_1
   \,|\, a_i \in H \rangle. \]
Then we have
\begin{theorem}[Satoh, \cite{S11}]\label{T-S11}
For any $k \geq 2$ and $n \geq k+2$,
\[ \mathrm{Coker}(\tau_{k, \Q}') \cong \mathcal{C}_n^{\Q}(k). \]
\end{theorem}

In our recent paper \cite{ES1}, which is a joint work with Naoya Enomoto,
we gave a combinatorial description of the $\mathrm{GL}(n,\Q)$-irreducible decomposition of $\mathcal{C}_n^{\Q}(k)$.
We remark that, as a $\mathrm{GL}(n,\Q)$-module, $\mathcal{C}_n^{\Q}(k)$ is isomorphic to the invariant part
$a_n(k):=(H_{\Q}^{\otimes k})^{\mathrm{Cyc}_k}$ of $H_{\Q}^{\otimes k}$
by the action of the cyclic group $\mathrm{Cyc}_k$ of order $k$. Namely, the cokernel $\mathrm{Coker}(\tau_{k, \Q}')$ is isomorphic to Kontsevich's
$a_n(k)$ as a $\mathrm{GL}(n,\Q)$-module. We also remark that in our notation $a_n(k)$ is considered for any $n \geq 2$ in contrast to Kontsevich's notation
for even $n=2g$. (See \cite{Kon1} and \cite{Kon2}.)

\vspace{0.5em}

Now we give tables of the irreducible decompositions of $\mathcal{C}_n^\Q(k)$ and $\mathrm{Im}(\tau_{k,\Q}')$ for $1 \leq k \leq 7$.

\vspace{0.5em}

{\small
\begin{center}
{\renewcommand{\arraystretch}{1.3}
\begin{tabular}{|c|l|l|} \hline
  $k$  & \hspace{3em} $\mathcal{C}_n^\Q(k)=\mathrm{Coker}(\tau_{k,\Q}')$, \,\,\, $n \geq k+2$        &                         \\ \hline
  $1$  & $0$                                                 & Andreadakis \cite{And}  \\ \hline
  $2$  & $(2)$                                           & Pettet \cite{Pet}       \\ \hline
  $3$  & $(3) \oplus (1^3)$                          & Satoh \cite{S03}       \\ \hline
  $4$  & $(4) \oplus (2,2) \oplus (2,1^2)$  & Satoh \cite{S09} \\ \hline
  $5$  & $(5) \oplus (3,2) \oplus 2(3,1^2) \oplus (2^2,1) \oplus (1^5)$ & \\ \hline
  $6$  & $(6) \oplus 2(4,2) \oplus 2(4,1^2) \oplus (3^2) \oplus 2(3,2,1) \oplus (3,1^3) \oplus 2(2^3)$ & \\
& \hspace{1.1em} $\oplus (2^2,1^2) \oplus (2,1^4)$ & \\ \hline
  $7$  & $(7) \oplus 2(5,2) \oplus 3(5,1^2) \oplus 2(4,3) \oplus 5(4,2,1) \oplus 2(4,1^3) \oplus 3(3^2,1)$ & \\
       & \hspace{1.1em} $\oplus 3(3,2^2)\oplus 5(3,2,1^2) \oplus 3(3,1^4) \oplus 2(2^3,1) \oplus 2(2^2,1^3) \oplus (1^7)$ & \\ \hline
\end{tabular}}
\end{center}
}

In the table above, $(\lambda)$ denotes the irreducible polynomial representation $L^{(\lambda)}$ of $\mathrm{GL}(n,\Q)$ corresponding to a Young tableau $\lambda$. (For the notation, see \cite{ES1} for details.)

{\small
\begin{align*}
{\renewcommand{\arraystretch}{1.2}
\begin{array}{|c|c|c|}
\hline
k & \text{polynomial part of} \ \mathrm{Im}(\tau_{k,\Q}') & \text{non-polynomial part of} \ \mathrm{Im}(\tau_{k,\Q}') \\
\hline
1 & (1) & (1,1) \\
\hline
2 & (1^2) & (2,1) \\
\hline
3 & 2(2,1) & (3,1) \oplus (2,1^2) \\
\hline
4 & 3(3,1) \oplus (2^2) \oplus 2(2,1^2) \oplus (1^4) & (4,1) \oplus (3,2) \oplus (3,1^2) \\
& & \oplus (2^2,1) \oplus (2,1^3)\\
\hline
5 & 4(4,1) \oplus 4(3,2) \oplus 4(3,1^2) & (5,1) \oplus (4,2) \oplus 2(4,1^2) \oplus (3^2) \\
 & \oplus 4(2^2,1) \oplus 4(2,1^3)& 3(3,2,1) \oplus (3,1^3) \oplus 2(2^2,1^2) \oplus (2,1^4)\\
\hline
6 & 5(5,1) \oplus 7(4,2) \oplus 8(4,1^2) \oplus 4(3^2) & (6,1) \oplus 2(5,2) \oplus 2(5,1^2) \oplus 2(4,3) \\
& \oplus 14(3,2,1) \oplus 9(3,1^3) \oplus 3(2^3)  & \oplus 5(4,2,1) \oplus 3(4,1^3) \oplus 3(3^2,1) \\
& \oplus 8(2^2,1^2) \oplus 4(2,1^4) \oplus (1^6) & \oplus 3(3,2^2) \oplus 5(3,2,1^2) \oplus 2(3,1^4) \\
& & \oplus 2(2^3,1) \oplus 2(2^2,1^3) \oplus (2,1^5)\\
\hline
7 & 6(6,1) \oplus 12(5,2) \oplus 12(5,1^2) \oplus 12(4,3) & (7,1) \oplus 2(6,2) \oplus 3(6,1^2) \oplus 4(5,3) \\
&  \oplus 30(4,2,1) \oplus 18(4,1^3) \oplus 18(3^2,1)  & \oplus 8(5,2,1) \oplus 4(5,1^3) \oplus (4^2) \oplus 9(4,3,1) \\
& \oplus 18(3,2^3) \oplus 30(3,2,1^2) \oplus 12(3,1^4) & \oplus 6(4,2^2) \oplus 12(4,2,1^2) \oplus 4(4,1^4) \oplus 6(3^2,2) \\
& \oplus 12(2^3,1) \oplus 12(2^2,1^3) \oplus 6(2,1^5) & \oplus 9(3,2^2,1) \oplus 8(3,2,1^3) \oplus 3(3,1^5) \oplus (2^4)\\
& & \oplus 4(2^3,1^2) \oplus 2(2^2,1^4) \oplus (2,1^6) \\
\hline
\end{array}}
\end{align*}
}

In the polynomial part of the table above, $(\lambda)$ denotes the irreducible polynomial representation $L^{(\lambda)}$,
and in the non-polynomial part, $(\mu)$ denotes the irreducible non-polynomial representation $L^{\{\mu;(1)\}}$
of $\mathrm{GL}(n,\Q)$ corresponding to Young tableaux $\lambda$ and $\mu$. (For the notation, see \cite{ES1} for details.)

\vspace{1em}

By Corollary {\rmfamily \ref{C-Mor}}, the irreducible module $S^k H_{\Q} = (k)$ appears in $\mathrm{Coker}(\tau_{k,\Q}')$ for each $k \geq 2$.
On the other hand, in \cite{S03} we showed that the irreducible module $\Lambda^k H_{\Q} = (1^k)$ appears in $\mathrm{Coker}(\tau_{k,\Q}')$ for each odd $k \geq 3$.
In \cite{ES1}, we showed that their multiplicities of them in $\mathrm{Coker}(\tau_{k,\Q}') \cong \mathcal{C}_n^\Q(k)$ are one. Namely,
\begin{theorem}[Enomoto and Satoh, \cite{ES1}]\label{T-ES1}
For $k \geq 2$ and $n \geq k+2$,
\begin{enumerate}
\item $[\mathrm{Coker}(\tau_{k,\Q}'): S^k H_{\Q}]=1$.
\item $[\mathrm{Coker}(\tau_{k,\Q}'): \Lambda^k H_{\Q}]=1$ if $k$ is odd.
\end{enumerate}
\end{theorem}

From this theorem, we see that the multiplicity of the Morita obstruction $S^k H_{\Q}$ in $\mathrm{Coker}(\tau_{k,\Q})$ is also one.
Hence, the Morita obstruction is essentially unique in the cokernel of the rational Johnson homomorphism $\tau_{k,\Q}$.

\section{Free metabelian group case}\label{S-FreM}

In this section, we consider a metabelian analogue of Section {\rmfamily \ref{S-Fre}}. The main purpose of this section is to give the irreducible
decomposition of the images and the cokernels of the Johnson homomorphisms of the automorphism groups of free metabelian groups.

\vspace{0.5em}

The quotient group of the free group $F_n$ by its second derived group $[[F_n,F_n],[F_n,F_n]]$ is called a free metabelian group of rank $n$,
denoted by $F_n^M$.
Then, $(F_n^M)^{\mathrm{ab}}=H$, and hence $\mathrm{Aut}\,(F_n^M)^{\mathrm{ab}}=\mathrm{Aut}(H)=\mathrm{GL}(n,\Z)$.
Since the surjective homomorphism $\mathrm{Aut}\,F_n \rightarrow \mathrm{GL}(n,\Z)$ factors through
$\mathrm{Aut}\,F_n^M$, a homomorphism $\mathrm{Aut}\,F_n^M \rightarrow \mathrm{GL}(n,\Z)$ induced from the abelianization of $F_n^M$ is surjective.
Thus, we also identify $\mathrm{Aut}\,F_n^M/\mathrm{IA}(F_n^M)$ with $\mathrm{GL}(n,\Z)$.

\vspace{0.5em}

In the following, for simplicity, we omit the capital \lq\lq $F$" from a symbol attached with $F_n^M$ if there is no confusion.
For example, we write $\Gamma_n^M(k)$, $\mathrm{gr}^k(\mathcal{L}_n^M)$, $\mathrm{IA}_n$, $\mathcal{A}_n^M(k)$ and $\mathrm{gr}^k(\mathcal{A}_n^M)$
for $\Gamma_{F_n^M}(k)$, $\mathrm{gr}^k(\mathcal{L}_{F_n^M})$, $\mathrm{IA}(F_n^M)$ and $\mathrm{gr}^k(\mathcal{A}_{F_n^M})$ respectively.

\subsection{Free metabelian Lie algebras and their derivations}\label{Ss-FMLA}

By a remarkable work of Chen \cite{Che}, it is known that each of $\mathrm{gr}^k(\mathcal{L}_n^M)$ is a free abelian group, and that
the graded Lie algebra $\mathrm{gr}(\mathcal{L}_n^M)$ is isomorphic to the free metabelian Lie algebra generated by $H$.
Chen \cite{Che} also gave the rank of $\mathrm{gr}^k(\mathcal{L}_n^M)$. That is,
\begin{equation}\label{ex-chen}
 \mathrm{rank}_{\Z}(\mathrm{gr}^k(\mathcal{L}_n^M))=(k-1) \binom{n+k-2}{k}.
\end{equation}

\vspace{0.5em}

Consider the derivation algebra $\mathrm{Der}(\mathrm{gr}(\mathcal{L}_n^M))$ of the free Lie algebra $\mathrm{gr}(\mathcal{L}_n^M)$.
By the universality of the free metabelian Lie algebra $\mathrm{gr}(\mathcal{L}_n^M)$, the embedding
\[ \mathrm{Der}(\mathrm{gr}(\mathcal{L}_n^M))(k) \hookrightarrow H^* {\otimes}_{\Z} \mathrm{gr}^{k+1}(\mathcal{L}_n^M) \]
as mentioned above is surjective.
Namely, the degree $k$-part $\mathrm{Der}(\mathrm{gr}(\mathcal{L}_n^M))(k)$ is considered as
\[ \mathrm{Hom}_{\Z}(H,\mathrm{gr}^{k+1}(\mathcal{L}_n^M)) = H^* {\otimes}_{\Z} \mathrm{gr}^{k+1}(\mathcal{L}_n^M) \]
for each $k \geq 1$.

\subsection{The IA-automorphism group of $F_n^M$}\label{Ss-IAFM}

Consider a homomorphism $\mathrm{Aut}\,F_n \rightarrow \mathrm{Aut}\,F_n^M$ induced from the action of
$\mathrm{Aut}\,F_n$ on $F_n^M$. Restricting it to $\mathrm{IA}_n$,
we obtain a homomorphism $\nu : \mathrm{IA}_n \rightarrow \mathrm{IA}_n^M$.
Bachmuth and Mochizuki \cite{BM1} showed that $\nu$ is not surjective for $n =3$, and that $\mathrm{IA}_3^M$ is not finitely generated.
They also showed in \cite{BM2} that $\nu$ is surjective for $n \geq 4$.
Hence for $n \geq 4$, $\mathrm{IA}_n^M$ is finitely generated by the images of the Magnus generators of $\mathrm{IA}_n$.
It is, however, not known whether $\mathrm{IA}_n^M$ is finitely presentable or not for $n \geq 4$.

\subsection{Johnson homomorphisms of $\mathrm{Aut}\,F_n^M$}\label{Ss-JohnM}

In the following, we consider the case where $n \geq 4$ unless otherwise noted.

\vspace{0.5em}

To begin with, let $\mathcal{K}_n$ be the kernel of the homomorphism $\nu : \mathrm{IA}_n \rightarrow \mathrm{IA}_n^M$.
Then by the definition of the Johnson filtration of $\mathrm{Aut}\,F_n$, we see that $\mathcal{K}_n \subset \mathcal{A}_n(2)=[\mathrm{IA}_n, \mathrm{IA}_n]$.
Hence, we obtain
\[ (\mathrm{IA}_n^M)^{\mathrm{ab}} \cong \mathrm{IA}_n^{\mathrm{ab}} \cong H^* \otimes_{\Z} \Lambda^2 H. \]
Furthermore, the first Johnson homomorphism $\tau_1$ of $\mathrm{Aut}\,F_n^M$ is just the abelianization of $\mathrm{IA}_n^M$.

\vspace{0.5em}

In \cite{S06}, by showing that the Morita trace $\mathrm{Tr}_{[k]} : H^* {\otimes}_{\Z} \mathrm{gr}^{k+1}(\mathcal{L}_n) \rightarrow S^k H$
factors through a homomorphism $\mathrm{gr}^k(\mathcal{A}_n) \rightarrow \mathrm{gr}^k(\mathcal{A}_n^M)$ induced from the inclusion,
we proved 
\begin{theorem}[Satoh, \cite{S06}]
For any $n \geq 4$ and $k \geq 2$,
\[ 0 \rightarrow {\mathrm{gr}}^k (\mathcal{A}_n^M) \xrightarrow{\tau_k} H^* \otimes_{\Z} \mathrm{gr}^{k+1}(\mathcal{L}_n^M)
     \xrightarrow{\mathrm{Tr}_{[k]}} S^k H \rightarrow 0 \]
is a $\mathrm{GL}(n,\Z)$-equivariant exact sequence.
\end{theorem}

Moreover, recently, we gave the irreducible decomposition of $\mathrm{gr}_{\Q}^k(\mathcal{A}_n^M)$ as follows:
\begin{theorem}[Enomoto and Satoh, \cite{ES1}]
For $k \geq 2$,
\[ \mathrm{gr}_{\Q}^k(\mathcal{A}_n^M) \cong L^{\{(k,1),(1)\}} \oplus L^{\{(k-1,1),0\}}. \]
\end{theorem}

\subsection{The abelianization of $\mathrm{Der}(\mathrm{gr}(\mathcal{L}_n^M))$}

Here we show an application of Morita's trace maps $\mathrm{Tr}_{[k]}$.
In \cite{ES1}, we determined the abelianization of the derivation algebra
$\mathrm{Der}^+(\mathrm{gr}(\mathcal{L}_n^M))$ as a Lie algebra using Morita's trace maps. That is,
\begin{theorem}[Enomoto and Satoh, \cite{ES1}]
For $n \geq 4$, we have
\[ (\mathrm{Der}^+(\mathcal{L}_n^M))^{\mathrm{ab}} \cong (H^* \otimes_{\Z} \Lambda^2 H ) \oplus \bigoplus_{k \geq 2} S^k H. \]
\end{theorem}

More precisely, we showed that this isomorphism is given by the degree one part of $\mathrm{Der}^+(\mathcal{L}_n^M)$
and Morita's trace maps $\mathrm{Tr}_{[k]}$ for $k \geq 2$.
Namely, consider a Lie algebra homomorphism
\[ \Theta^M := \mathrm{id}_1 \oplus \bigoplus_{k \geq 2} \mathrm{Tr}_{[k]} : \mathrm{Der}^+(\mathcal{L}_n^M) \rightarrow (H^* \otimes_{\Z} \Lambda^2 H) \oplus
   \bigoplus_{k \geq 2} S^k H. \]
Then we obtained an exact sequence
\[ 0 \rightarrow \bigoplus_{k \geq 2} \mathrm{gr}^k(\mathcal{A}_n^M) \xrightarrow{\oplus_{k \geq 2} \tau_k}
      \mathrm{Der}^+(\mathcal{L}_n^M) \xrightarrow{\Theta^M} (H^* \otimes_{\Z} \Lambda^2 H ) \oplus \bigoplus_{k \geq 2} S^k H \rightarrow 0. \]
of Lie algebras.

\vspace{0.5em}

We should remark that this result is the Chen Lie algebra version of Morita's conjecture for the free Lie algebra $\mathrm{gr}(\mathcal{L}_n)$:
\begin{conjecture}[Morita]
For $n \geq 3$,
\[ (\mathrm{Der}^+(\mathcal{L}_n))^{\mathrm{ab}} \cong (H^* \otimes_{\Z} \Lambda^2 H ) \oplus \bigoplus_{k \geq 2} S^k H \]
\end{conjecture}
Morita \cite{Mo3} showed that this conjecture is true up to degree $n(n-1)$, based on a work of Kassabov \cite{Kas}.

\subsection{On the second cohomology of $\mathrm{IA}_n^M$}\label{Ss-Cup}

Due to Pettet \cite{Pet}, the image of the cup product
\[ \cup : \Lambda^2 H^1(\mathrm{IA}_n,\Q) \rightarrow H^2(\mathrm{IA}_n,\Q) \]
is completely determined as mentioned above.
In \cite{S06}, we study a metabelian analogue of her work. 
Let
\[ \cup^M : \Lambda^2 H^1(\mathrm{IA}_n^M,\Q) \rightarrow H^2(\mathrm{IA}_n^M,\Q) \]
be the cup product of $H^1(\mathrm{IA}_n^M,\Q)$. Then, we proved
\begin{theorem}[Satoh, \cite{S06}]
For $n \geq 4$, the natural homomorphism $\nu : \mathrm{IA}_n \rightarrow \mathrm{IA}_n^M$ induces an isomorphism
\[ \nu^* : \mathrm{Im}(\cup^M) \rightarrow \mathrm{Im}(\cup). \]
\end{theorem}

Observing the five term exact sequence of a group extension
\[ 1 \rightarrow \mathcal{K}_n \rightarrow \mathrm{IA}_n \rightarrow \mathrm{IA}_n^M \rightarrow 1, \]
we have
\[\begin{split}
   0 & \rightarrow H^1(\mathrm{IA}_n^M,\Q) \rightarrow H^1(\mathrm{IA}_n,\Q) \\
     & \hspace{2em} \rightarrow H^1(\mathcal{K}_n,\Q)^{\mathrm{IA}_n}
          \rightarrow H^2(\mathrm{IA}_n^M,\Q) \rightarrow H^2(\mathrm{IA}_n,\Q).
  \end{split}\]
Since $H^1(\mathrm{IA}_n^M,\Q) \rightarrow H^1(\mathrm{IA}_n,\Q)$ is isomorphism by $(\mathrm{IA}_n^M)^{\mathrm{ab}} \cong \mathrm{IA}_n^{\mathrm{ab}}$
as mentioned above, we obtain
\[ 0 \rightarrow H^1(\mathcal{K}_n,\Q)^{\mathrm{IA}_n}
          \rightarrow H^2(\mathrm{IA}_n^M,\Q) \rightarrow H^2(\mathrm{IA}_n,\Q). \]
In \cite{S06}, we showed using the third Johnson homomorphism $\tau_3$ of $\mathrm{Aut}\,F_n^M$ that $H^1(\mathcal{K}_n,\Q)^{\mathrm{IA}_n}$ is non-trivial.
Hence we see
\begin{theorem}[Satoh, \cite{S06}]
For $n \geq 4$, $\cup^M$ is not surjective.
\end{theorem}

Recently, we \cite{S14} showed that $\mathcal{K}_n$ is not finitely generated. This will be mentioned more precisely
in Subsection {\rmfamily \ref{Ss-KMag}} later.

\section{Braid groups and Mapping class groups}\label{S-BM}

In this section, we recall well known facts of
two important subgroups of the automorphism group of a free group.
The first one is the braid group, and the other is the mapping class group of a surface.
Needless to say, these groups have been studied by a large number of authors with a long history.

\subsection{Braid groups}

Let $B_n$ be the {\it braid group} \index{braid group} of $n$ strands. Artin \cite{Art} gave the first finite presentation of
$B_n$ with generators $\sigma_i$ for $1 \leq i \leq n-1$ subject to relations:
\begin{itemize}
\item $\sigma_i \sigma_{i+1} \sigma_i = \sigma_{i+1} \sigma_{i} \sigma_{i+1}$ for $1 \leq i \leq n-2$, 
\item $[\sigma_i, \sigma_j]=1$ for $|i-j| \geq 2$.
\end{itemize}
It is known that the braid group $B_n$ is isomorphic to the mapping class group of the $n$-punctured disk
$D^2 \setminus \{ p_1, \ldots p_n \}$. (See Theorem {\rmfamily 1.10} in \cite{Bir}.)
Then the action of $B_n$ on the fundamental group $\pi_1(D^2 \setminus \{ p_1, \ldots p_n \}) \cong F_n$
induces a faithful representation of $B_n$ as the automorphisms of $F_n$. This representation, denoted by
$\psi : B_n \hookrightarrow \mathrm{Aut}\,F_{n}$, is given by
\[ \sigma_i : x_i \mapsto x_i x_{i+1} x_i^{-1}, \hspace{1em} x_{i+1} \mapsto x_{i}, \hspace{1em} x_j \mapsto x_j \]
for $j \neq i$, $i+1$.
Furthermore, its image is completely determined as
\begin{theorem}[See Theorem {\rmfamily 1.9} in \cite{Bir}]\label{T-braid}
\[\begin{split}
  \{ \sigma \in \mathrm{Aut}\,F_{n} \,\, | \,\, x_i^{\sigma}=A_i \, x_{\mu(i)} \, A_i^{-1}, \,\,
     \mu \in \mathfrak{S}_n, \, A_i \in F_n, \,\, \omega^{\sigma}=\omega \}
  \end{split}\]
where $\mathfrak{S}_n$ is the symmetric group of degree $n$,
and $\omega=x_1x_2 \cdots x_n \in F_n$. Namely, $\omega$ is a homotopy class of a simple closed curve
on $D^2 \setminus \{ p_1, \ldots p_n \}$ parallel to the boundary.
\end{theorem}

Through this embedding, we consider $B_n$ as a subgroup of $\mathrm{Aut}\,F_{n}$.
Using this embedding, Artin solved the word problem for the braid group $B_n$. (See \cite{Bir} for details.)
A subgroup $B_n \cap \mathrm{IA}_n$ is called the {\it pure braid group} \index{pure braid group}, and is denoted by $P_n$.

\subsection{Mapping class groups of surfaces}

Next, let us consider the mapping class groups of surfaces.
(See also for example \cite{Mo6} by Morita for basic material for the mapping class groups.)
For any integer $g \geq 1$, let $\Sigma_{g,1}$ be the compact oriented surface of genus $g$ with one boundary component.
The fundamental group $\pi_1(\Sigma_{g,1})$ of $\Sigma_{g,1}$ is isomorphic to the free group $F_{2g}$.
Throughout the paper, we fix isomorphisms $\pi_1(\Sigma_{g,1}) \cong F_{2g}$ and
$H_1(\Sigma_{g,1},\Z) \cong H$.

\vspace{0.5em}

We denote by $\mathcal{M}_{g,1}$ the {\it mapping class group} \index{mapping class group} of ${\Sigma}_{g,1}$.
Namely, $\mathcal{M}_{g,1}$ is the group of isotopy classes of
orientation preserving diffeomorphisms of ${\Sigma}_{g,1}$ which fix the boundary pointwise.
Then the action of $\mathcal{M}_{g,1}$ on $\pi_1(\Sigma_{g,1})=F_{2g}$ induces a natural homomorphism
\[ \phi : \mathcal{M}_{g,1} \hookrightarrow \mathrm{Aut}\,F_{2g}. \]
By classical works due to Dehn and Nielsen, we have
\begin{theorem}[Dehn and Nielsen]\label{T-DN}
For any $g \geq 1$, $\phi$ is injective. Furthermore, its image is given by
\[ \mathrm{Im}(\phi) = \{ \sigma \in \mathrm{Aut}\,F_{2g} \,\, | \,\, \zeta^{\sigma} = \zeta \} \]
where $\zeta \in F_{2g}$ is the homotopy class of a simple closed curve on $\Sigma_{g,1}$ which is parallel to the boundary.
\end{theorem}

By this embedding, we consider $\mathcal{M}_{g,1}$ as a subgroup of $\mathrm{Aut}\,F_{2g}$.
The image of $\mathcal{M}_{g,1}$ by $\rho : \mathrm{Aut}\,F_{2g} \rightarrow \mathrm{GL}(2g,\Z)$ is the symplectic group
\[ \mathrm{Sp}(2g,\Z):=\{X \in \mathrm{GL}(2g,\Z) \ | \ {}^{t}X J X = J \} \,\,\, \mathrm{for} \,\,\, J= \left(
\begin{array}{cc}
0 & E_g \\
-E_g & 0
\end{array}
\right) \]
where $E_g$ is the identity matrix of degree $g$.

\vspace{0.5em}

Set $\mathcal{I}_{g,1} := \mathcal{M}_{g,1} \cap \mathrm{IA}_n$. Namely, $\mathcal{I}_{g,1}$ is
a subgroup of $\mathcal{M}_{g,1}$ consisting of mapping classes which act on the first homology $H_1(\Sigma_{g,1},\Z)$ of
$\Sigma_{g,1}$ trivially. The group $\mathcal{I}_{g,1}$ is called the {\it Torelli subgroup} of $\mathcal{M}_{g,1}$,
or the {\it Torelli group} \index{Torelli group} of $\Sigma_{g,1}$. Then we have a commutative diagram including three group extensions:
\[\begin{CD}
  1 @>>> P_n @>>> B_n @>>> \mathfrak{S}_n @>>> 1 \\
  @.    @VVV   @VV{\psi}V       @VVV             @. \\
  1 @>>> \mathrm{IA}_n @>>> \mathrm{Aut}\,F_n @>>> \mathrm{GL}(n,\Z) @>>> 1 \\
  @.    @A{n=2g}AA   @A{n=2g}A{\phi}A       @A{n=2g}AA             @. \\
  1 @>>> \mathcal{I}_{g,1} @>>> \mathcal{M}_{g,1} @>>> \mathrm{Sp}(2g,\Z) @>>> 1
  \end{CD}\]

\section{Magnus representations}\label{S-Mag}

In this section we consider the Magnus representation of $\mathrm{Aut}\,F_n$ and $\mathrm{Aut}\,F_n^M$.
The main goal is to show a relation between the Johnson homomorphism and the Magnus representation due to Morita \cite{Mo1},
and the fact that the kernel of the Magnus representation
is not finitely generated. For basic material for the Magnus representation, see for example \cite{Bir} or \cite{Sak}.

\subsection{Fox derivations}

To begin with, we recall some properties of the Fox derivations of $F_n$. See also \cite{Sak}.

\begin{definition}
For each $1 \leq i \leq n$, let
\[ \frac{\partial}{\partial x_i} : \Z[F_n] \rightarrow \Z[F_n] \]
be a map defined by
\[ \frac{\partial}{\partial x_i}(w) = 
                            \sum_{j=1}^{r} \epsilon_j \delta_{\mu_j,i} x_{\mu_1}^{\epsilon_1} \cdots
      x_{\mu_j}^{\frac{1}{2}(\epsilon_j-1)} \in \Z[F_n] \]
for a reduced word $w=x_{\mu_1}^{\epsilon_1} \cdots x_{\mu_r}^{\epsilon_r} \in F_n$, $\epsilon_j=\pm1$.
The maps $\frac{\partial}{\partial x_i}$ are called the {\it Fox derivations} \index{Fox derivation} of $F_n$.
\end{definition}

\begin{definition}
For a group $G$, let $\Z[G]$ be the integral group ring of $G$ over $\Z$. We denote the augmentation map
by $\epsilon : \Z[G] \rightarrow \Z$.
The kernel $I_G$ of $\epsilon$ is called the augmentation ideal of $\Z[G]$. Then the $k$-times powers $I_G^k:=I_G \times \cdots \times I_G$
of $I_G$ for all $k \geq 1$ provide
a descending filtration of $\Z[G]$, and the direct sum
\[ \mathrm{gr}(\Z[G]) := \bigoplus_{k \geq 1} \, I_G^k/I_G^{k+1} \]
naturally has a graded algebra structure induced from the multiplication of $\Z[G]$.
It is called the {\it graded algebra associated to the group ring} $\Z[G]$. \index{graded algebra associated to a group ring}
\end{definition}

For $G=F_n$, the free group of rank $n$,
it is classically known (by work of Magnus \cite{MKS}) that each graded quotient $I_{F_n}^k/I_{F_n}^{k+1}$ is a free abelian group with basis
\[ \{ (x_{i_1}-1)(x_{i_2}-1) \cdots (x_{i_k}-1) \,\,|\,\, 1 \leq i_j \leq n \}, \]
and a map $I_{F_n}^k/I_{F_n}^{k+1} \rightarrow H^{\otimes k}$ defined by
\[ (x_{i_1}-1)(x_{i_2}-1) \cdots (x_{i_k}-1) \mapsto x_{i_1} \otimes x_{i_2} \otimes \cdots \otimes x_{i_k} \]
induces an associative algebra isomorphism from $\mathrm{gr}(\Z[F_n])$ to the tensor algebra $T(H) := \bigoplus_{k \geq 1} H^{\otimes k}$.
In the following, we identify $I_{F_n}^k/I_{F_n}^{k+1}$ with $H^{\otimes k}$ via this isomorphism.

\vspace{0.5em}

In $\mathrm{gr}(\Z[F_n])$, the free Lie algebra $\mathrm{gr}(\mathcal{L}_n)$ is realized as follows.
\begin{theorem}[Fox, \cite{Fo1}] \quad
\begin{enumerate}
\item For any $k \geq 1$ and $x \in \Gamma_n(k)$, $x-1 \in I_{F_n}^k$. Furthermore, a map $\Gamma_n(k) \rightarrow I_{F_n}^k$ defined by $x \mapsto x-1$
induces an injective homomorphism
\[ \mathrm{gr}^k(\mathcal{L}_n) \hookrightarrow I_{F_n}^k/I_{F_n}^{k+1}. \]
\item For any $k \geq 1$ and $x \in \Gamma_n(k)$, $\frac{\partial x}{\partial x_i} \in I_{F_n}^{k-1}$ for any $1 \leq i \leq n$. Furthermore,
\[ x-1 = \sum_{i=1}^n \frac{\partial x}{\partial x_i} (x_i-1) \in \Z[F_n]. \]
\end{enumerate}
\end{theorem}

The formula Part (2) above is called the fundamental formula of the Fox derivations.

\vspace{0.5em}

On the other hand, for $G=H$, the free abelian group of rank $n$,
it is also known that each graded quotient $I_H^k/I_H^{k+1}$ is a free abelian group with basis
\[ \{ (x_{i_1}-1)(x_{i_2}-1) \cdots (x_{i_k}-1) \,\,|\,\, 1 \leq i_1 \leq i_2 \leq \cdots \leq i_k \leq n \}, \]
and the associated graded algebra $\mathrm{gr}(\Z[H])$ is isomorphic to the symmetric algebra
\[ S(H) := \bigoplus_{k \geq 1} S^k H \]
of $H$. (See \cite{Pas}, Chapter VIII, Proposition 6.7.) We also identify $I_H^k/I_H^{k+1}$ with $S^k H$.

\vspace{0.5em}

Let $\mathfrak{a} : F_n \rightarrow H$ be the abelianization of $F_n$. 
Then a homomorphism $I_{F_n}^k/I_{F_n}^{k+1} \rightarrow I_H^k/I_H^{k+1}$ induced from the abelianization
$\mathfrak{a} : F_n \rightarrow H$ is just the natural projection $H^{\otimes k} \rightarrow S^k H$.

\subsection{Magnus representations of $\mathrm{Aut}\,F_n$}\label{Ss-MJ}

Let $\varphi : F_n \rightarrow G$ be a group homomorphism.
We also denote by $\varphi$ the ring homomorphism
$\Z[F_n] \rightarrow \Z[G]$ induced from $\varphi$.
For any matrix $A=(a_{ij}) \in \mathrm{GL}(n,\Z[F_n])$, let $A^{\varphi}$ be the matrix
$(a_{ij}^{\varphi}) \in \mathrm{GL}(n,\Z[G])$.

\begin{definition}
A map $r^{\varphi} : \mathrm{Aut}\,F_n \rightarrow \mathrm{GL}(n,\Z[G])$
defined by
\[ \sigma \mapsto \biggl{(} \frac{\partial {x_i}^{\sigma} }{\partial x_j} {\biggl{)}}^{\varphi} \]
for any $\sigma \in \mathrm{Aut}\,F_n$, is called the {\it Magnus representation} \index{Magnus representation}
of $\mathrm{Aut}\,F_n$ associated to $\varphi$.
If $\varphi= \mathrm{id} : F_n \rightarrow F_n$, we write $r$ for $r^{\mathrm{id}}$ for simplicity.
\end{definition}

For the abelianization $\mathfrak{a} : F_n \rightarrow H$,
the map $r^{\mathfrak{a}}$ is not a homomorphism but a crossed homomorphism (with respect to right action). Namely,
$r^{\mathfrak{a}}$ satisfies
\[ r^{\mathfrak{a}}(\sigma \tau) = (r^{\mathfrak{a}}(\sigma))^{\tau^*} \cdot r^{\mathfrak{a}}(\tau) \]
for any $\sigma, \tau \in \mathrm{Aut}\,F_n$.
Here $(r^{\mathfrak{a}}(\sigma))^{\tau^*}$ denotes the matrix obtained from $r^{\mathfrak{a}}(\sigma)$
by applying the automorphism
$\tau^* : \Z[H] \rightarrow \Z[H]$ induced from $\rho(\tau) \in \mathrm{Aut}(H)$ on each entry.
Hence by restricting $r^{\mathfrak{a}}$ to $\mathrm{IA}_n$, we obtain a group homomorphism
$\mathrm{IA}_n \rightarrow \mathrm{GL}(n,\Z[H])$, which is also denoted by $r^{\mathfrak{a}}$.

\vspace{0.5em}

Here we remark a relation between the Magnus representation and the Johnson homomorphism of $\mathrm{Aut}\,F_n$.
For any element $f \in H^* \otimes_{\Z} \mathrm{gr}^{k+1}(\mathcal{L}_n)$, set
\[ \begin{Vmatrix} f \end{Vmatrix}
    := \biggl{(} \frac{\partial {({x_i}^{f})} }{\partial x_j} {\biggl{)}} \in M(n, I_{F_n}^k/I_{F_n}^{k+1}) \]
where we consider any lift of the element
\[ {x_i}^{f} \in \mathrm{gr}^{k+1}(\mathcal{L}_n) = \Gamma_n(k+1)/\Gamma_n(k+2) \]
to $\Gamma_n(k+1)$. We can see that the equivalence class of $\frac{\partial {({x_i}^{f})} }{\partial x_j}$
in $I_{F_n}^k/I_{F_n}^{k+1}$ does not depend on the choice of the lift of ${x_i}^{f}$.
Then we have
\begin{theorem}[Morita, \cite{Mo1}]
The Magnus representation $r:\mathrm{Aut}\,F_n \rightarrow \mathrm{GL}(n,\Z[F_n])$ induces
a homomorphism
\[ r_k : \mathcal{A}_n(k) \rightarrow \mathrm{GL}(n,\Z[F_n]/I_{F_n}^{k+1}). \]
Moreover, it has the following form. For any $\sigma \in \mathcal{A}_n(k)$, we have
\[ r_k(\sigma)= E_n + \begin{Vmatrix} \tau_k(\sigma) \end{Vmatrix} \]
where $E_n$ is the identity matrix of degree $n$.
Namely, the representation $r_k$ is essentially equal to the $k$-th Johnson homomorphism.
\end{theorem}

\vspace{0.5em}

Finally, we remark Morita's original definition of the trace map $\mathrm{Tr}_{[k]}$.
For any element $f \in H^* \otimes_{\Z} \mathrm{gr}^{k+1}(\mathcal{L}_n)$, take the trace $\mathrm{Trace}(\begin{Vmatrix} f \end{Vmatrix}) \in I_{F_n}^k/I_{F_n}^{k+1}$
of the matrix $\begin{Vmatrix} f \end{Vmatrix}$.
Then we have
\[ \mathrm{Tr}_{[k]}(f) = (-1)^k (\mathrm{Trace}(\begin{Vmatrix} f \end{Vmatrix}))^{\mathfrak{a}} \in I_H^k/I_H^{k+1} \cong S^k H. \]
Since $\mathcal{A}_n(2)=\Gamma_{\mathrm{IA}_n}(2)$ by a result of Bachmuth \cite{Ba2}, we have
$\mathrm{det} \circ r^{\mathfrak{a}}(\sigma) = 1$ for any $\sigma \in \mathcal{A}_n(2)$. 
This shows that
\[ 1 = \mathrm{det} \circ r_k^{\mathfrak{a}}(\sigma) = 1 + \mathrm{Tr}_{[k]}(\tau_k(\sigma)) \in \Z[H]/I_H^{k+1} \]
for any $k \geq 2$ and $\sigma \in \mathcal{A}_n(k)$. Hence $\mathrm{Tr}_{[k]}(\tau_k(\sigma))=0 \in I_H^k/I_H^{k+1}$.
This means that the Morita trace $\mathrm{Tr}_{[k]}$ vanishes on the image of $\tau_k$.

\subsection{The kernel of the Magnus representation of $r^{\mathfrak{a}}$}\label{Ss-KMag}

Here we consider the kernel $\mathrm{Ker}(r^{\mathfrak{a}})$ of the Magnus representation $r^{\mathfrak{a}}$.
First, we recall a remarkable classical work of Bachmuth.
\begin{theorem}[Bachmuth, \cite{Ba1}]
For any $n \geq 2$, the Magnus representation $r^{\mathfrak{a}}$ factors through the natural homomorphism
$\nu : \mathrm{IA}_n \rightarrow \mathrm{IA}_n^M$. Moreover, the induced representation
\[ \mathrm{IA}_n^M \rightarrow \mathrm{GL}(n,\Z[H]) \]
of $\mathrm{IA}_n^M$ by $r^{\mathfrak{a}}$ is faithful.
\end{theorem}

This shows that $\mathrm{Ker}(r^{\mathfrak{a}}) = \mathcal{K}_n$, the kernel of $\nu : \mathrm{IA}_n \rightarrow \mathrm{IA}_n^M$.
Namely, the metabelianization of $F_n$ induces the injectivization of the Magnus representation $r^{\mathfrak{a}}$.

\vspace{0.5em}

Historically, the faithfulness of the Magnus representation restricted to two subgroups related to topology has been studied.
One is the pure braid group $P_n$, and the other is the Torelli group $\mathcal{I}_{g,1}$.
The restriction $r^{\mathfrak{a}}|_{P_n}$ of $r^{\mathfrak{a}}$ to $P_n$ is called the {\it Gassner representation}.
\index{Gassner representation} (cf. \cite{Sak}.)
For $n = 2$ and $3$, Magnus and Peluso \cite{MgP} showed that $r^{\mathfrak{a}}|_{P_n}$ is faithful.
It is, however, still an open problem to determine whether $r^{\mathfrak{a}}|_{P_n}$ is faithful or not for $n \geq 4$.

\vspace{0.5em}

On the other hand, it is known that $r^{\mathfrak{a}}|_{\mathcal{I}_{g,1}}$ is not faithful as follows:
\begin{theorem}[Suzuki, \cite{Su1}]
The Magnus representation $r^{\mathfrak{a}}|_{\mathcal{I}_{g,1}}$ restricted to the Torelli group $\mathcal{I}_{g,1}$ is not faithful
for $g \geq 2$.
\end{theorem}
In \cite{Su1}, Suzuki gave a non-trivial element in $\mathrm{Ker}(r^{\mathfrak{a}}|_{\mathcal{I}_{g,1}})$. Recently, Church and Farb gave the following remarkable
result for the infinite generation for $r^{\mathfrak{a}}|_{\mathcal{I}_{g,1}}$.
\begin{theorem}[Church and Farb, \cite{ChF}]
$H_1(\mathrm{Ker}(r^{\mathfrak{a}}|_{\mathcal{I}_{g,1}}),\Z)$ has infinite rank for $g \geq 2$.
\end{theorem}
This shows that the above $r^{\mathfrak{a}}|_{\mathcal{I}_{g,1}}$ is not finitely generated for $g \geq 2$.
They proved this theorem by using a result of Suzuki and a variant of the first Johnson homomorphism
of the mapping class group.

\vspace{0.5em}

Now, it is a natural question to ask how large is $\mathrm{Ker}(r^{\mathfrak{a}})$. Our answer is
\begin{theorem}[Satoh, \cite{S14}]
For any $n \geq 2$, $H_1(\mathrm{Ker}(r^{\mathfrak{a}}),\Z)$ has infinite rank.
\end{theorem}
To prove this, we consider some embeddings from $\mathrm{IA}_n$ into $\mathrm{IA}_m$ for various $m>n$, which
arise from the action of $\mathrm{IA}_n$ on finite-index subgroups of $F_n$.
Then using the first Johnson homomorphisms on $\mathrm{IA}_m$ which do not vanish on $\mathrm{Ker}(r^{\mathfrak{a}})$,
we can detect infinitely many independent elements in $H_1(\mathrm{Ker}(r^{\mathfrak{a}}),\Z)$. 
We remark that we can also show the above theorem by using the techniques of Church and Farb for even $n \geq 2$.

\section{Johnson homomorphisms of $\mathcal{M}_{g,1}$}\label{S-MCG}

In this section, we consider the Johnson homomorphisms of the mapping class groups of surfaces.
There is a broad range of results for the Johnson filtration and the Johnson homomorphisms of the mapping class groups.
In this chapter, we concentrate on the study of the cokernel of the Johnson homomorphisms.
By Dehn and Nielsen's classical work, we consider $\mathcal{M}_{g,1}$ as a subgroup
of $\mathrm{Aut}\,F_{2g}$ as above.

\begin{definition}
For each $k \geq 1$, set $\mathcal{M}_{g,1}(k) := \mathcal{M}_{g,1} \cap \mathcal{A}_{2g}(k)$.
Then we have a descending filtration
\[ \mathcal{I}_{g,1} = \mathcal{M}_{g,1}(1) \supset \mathcal{M}_{g,1}(2) \supset \cdots \supset \mathcal{M}_{g,1}(k) \supset \cdots \]
of the Torelli group $\mathcal{I}_{g,1}$. This filtration is called the {\it Johnson filtration} \index{Johnson filtration} of $\mathcal{M}_{g,1}$.
\end{definition}

Set $\mathrm{gr}^k (\mathcal{M}_{g,1}) := \mathcal{M}_{g,1}(k)/\mathcal{M}_{g,1}(k+1)$.
For each $k \geq 1$, the mapping class group $\mathcal{M}_{g,1}$ acts on $\mathrm{gr}^k (\mathcal{M}_{g,1})$ by conjugation.
This action induces that of $\mathrm{Sp}(2g,\Z)=\mathcal{M}_{g,1}/\mathcal{I}_{g,1}$ on it.
By an argument similar to that of $\mathrm{Aut}\,F_n$, the Johnson homomorphisms of $\mathcal{M}_{g,1}$ are defined as follows.
\begin{definition}
For $n=2g$ and $k \geq 1$, consider the restriction of
$\tilde{\tau}_k : \mathcal{A}_{2g}(k) \rightarrow \mathrm{Hom}_{\Z}(H, \mathrm{gr}^{k+1}(\mathcal{L}_{2g}))$ to $\mathcal{M}_{g,1}(k)$.
Then its kernel is just $\mathcal{M}_{g,1}(k+1)$. Hence we obtain an injective homomorphism
\[ \tau_k^{\mathcal{M}} : \mathrm{gr}^k (\mathcal{M}_{g,1}) \hookrightarrow \mathrm{Hom}_{\Z}(H, \mathrm{gr}^{k+1}(\mathcal{L}_{2g}))
       = H^* \otimes_{\Z} \mathrm{gr}^{k+1}(\mathcal{L}_{2g}). \]
The homomorphism $\tau_k^{\mathcal{M}}$ is $\mathrm{Sp}(2g,\Z)$-equivariant, and is called the $k$-th
{\it Johnson homomorphism} \index{Johnson homomorphism of the mapping class group} of $\mathcal{M}_{g,1}$.
\end{definition}
If we consider a $\mathrm{GL}(2g,\Z)$-module $H$ as an $\mathrm{Sp}(2g,\Z)$-module, then $H^* \cong H$ by Poincar\'{e} duality.
Hence, in the following, we canonically identify the target $H^* \otimes_{\Z} \mathrm{gr}^{k+1}(\mathcal{L}_{2g})$ of $\tau_k^{\mathcal{M}}$ with
$H \otimes_{\Z} \mathrm{gr}^{k+1}(\mathcal{L}_{2g})$.

\vspace{0.5em}

The Johnson filtration and the Johnson homomorphisms of $\mathcal{M}_{g,1}$ begun to be studied by D. Johnson \cite{Jo1} in the 1980s
who determined the abelianization of the Torelli group in \cite{Jo4} using the first Johnson homomorphism and the Birman-Craggs homomorphism.
In particular, he showed
\begin{theorem}[Johnson, \cite{Jo4}]
For $g \geq 3$, $\mathrm{Im}(\tau_1^{\mathcal{M}}) \cong \Lambda^3 H$ as an $\mathrm{Sp}(2g,\Z)$-module, and it gives the free part of $H_1(\mathcal{I}_{g,1},\Z)$.
\end{theorem}

\vspace{0.5em}

Now, let us recall the fact that the image of $\tau_k^{\mathcal{M}}$ is contained in a certain $\mathrm{Sp}(2g,\Z)$-submodule of
$H \otimes_{\Z} \mathrm{gr}^{k+1}(\mathcal{L}_{2g})$,
by a result of Morita \cite{Mo1}.
In general, for any $n \geq 1$,
let $H \otimes_{\Z} \mathrm{gr}^{k+1}(\mathcal{L}_{n}) \rightarrow \mathrm{gr}^{k+2}(\mathcal{L}_{n})$ be a $\GL(n,\Z)$-equivariant homomorphism defined by
\[ a \otimes X \mapsto [a,X], \hspace{1em} \mathrm{for} \hspace{1em} a \in H, \,\,\, X \in \mathrm{gr}^{k+1}(\mathcal{L}_{n}). \]
For $n=2g$, we denote by $\mathfrak{h}_{g,1}(k)$ the kernel of this homomorphism:
\[ \mathfrak{h}_{g,1}(k) := \mathrm{Ker}(H \otimes_{\Z} \mathrm{gr}^{k+1}(\mathcal{L}_{2g}) \rightarrow \mathrm{gr}^{k+2}(\mathcal{L}_{2g})). \]
Then Morita showed
\begin{theorem}[Morita, \cite{Mo1}]
For $k \geq 2$,
the image $\mathrm{Im}(\tau_k^{\mathcal{M}})$ of $\tau_k^{\mathcal{M}}$ is contained in $\mathfrak{h}_{g,1}(k)$.
\end{theorem}
Therefore, to determine how different is $\mathrm{Im}(\tau_k^{\mathcal{M}})$ from $\mathfrak{h}_{g,1}(k)$ is one of the most basic problems. 
Throughout theis chapter, the cokernel $\mathrm{Coker}(\tau_k^{\mathcal{M}})$ of $\tau_k^{\mathcal{M}}$ always means
the quotient $\mathrm{Sp}(2g,\Z)$-module $\mathfrak{h}_{g,1}(k)/\mathrm{Im}(\tau_k^{\mathcal{M}})$.

\vspace{0.5em}

In order to take advantage of representation theory, we consider the rationalization of modules. Namely, we study the $\mathrm{Sp}(2g,\Q)$-module
structures of $\mathrm{Im}(\tau_{k,\Q}^{\mathcal{M}})$ and $\mathrm{Coker}(\tau_{k,\Q}^{\mathcal{M}})$.
So far, the $\mathrm{Sp}$-irreducible decompositions of
$\mathrm{Im}(\tau_{k,\Q}^{\mathcal{M}})$ and  $\mathrm{Coker}(\tau_{k,\Q}^{\mathcal{M}})$ are determined for $1 \leq k \leq 4$ as follows.

{\small
\begin{center}
{\renewcommand{\arraystretch}{1.3}
\begin{tabular}{|c|l|l|l|} \hline
  $k$ & $\mathrm{Im}(\tau_{k,\Q}^{\mathcal{M}})$            & $\mathrm{Coker}(\tau_{k,\Q}^{\mathcal{M}})$   &                         \\ \hline
  $1$ & $[1^3] \oplus [1]$                                  & $0$                                  & Johnson \cite{Jo1}      \\ \hline
  $2$ & $[2^2] \oplus [1^2] \oplus [0]$                     & $0$                                  & Morita \cite{Mo0}, Hain \cite{Hai} \\ \hline
  $3$ & $[3,1^2] \oplus [2,1]$                              & $[3]$                                & Asada-Nakamura \cite{AN1}, Hain \cite{Hai}  \\ \hline
  $4$ & $[4,2] \oplus [3,1^3] \oplus [2^3]$ & $[2,1^2] \oplus [2]$  &   Morita \cite{Mo2}     \\
    & $\oplus 2 [3,1] \oplus [2,1^2] \oplus 2[2]$  &         &   \\ \hline
\end{tabular}}
\end{center}
}

\vspace{0.5em}

Morita \cite{Mo1} showed that the symmetric tensor product $S^k H_{\Q}$ appears in the $\mathrm{Sp}$-irreducible decomposition of
$\mathrm{Coker}(\tau_{k,\Q}^{\mathcal{M}})$
for odd $k \geq 3$ using the Morita trace map mentioned above. Moreover, Hiroaki Nakamura, in unpublished work, showed that
its multiplicity is one.
In general, however, to determine the cokernel of $\tau_k^{\mathcal{M}}$ is quite a difficult problem.

\vspace{1em}

Here, we recall a remarkable result of Hain.
As an $\mathrm{Sp}(2g,\Z)$-module, we consider $\mathfrak{h}_{g,1}(k)$ as a submodule of the degree $k$ part $\mathrm{Der}(\mathrm{gr}(\mathcal{L}_n))(k)$
of the derivation algebra of $\mathrm{gr}(\mathcal{L}_n)$. Moreover, the graded sum
\[ \mathfrak{h}_{g,1} := \bigoplus_{k \geq 1} \mathfrak{h}_{g,1}(k) \]
naturally has a Lie subalgebra structure of $\mathrm{Der}^+(\mathrm{gr}(\mathcal{L}_n))$. Therefore we obtain a graded Lie algebra homomorphism
\[ \tau^{\mathcal{M}} := \bigoplus_{k \geq 1} \tau_k^{\mathcal{M}} : {\mathrm{gr}}(\mathcal{M}_{g,1}) \rightarrow \mathfrak{h}_{g,1}, \]
which is called the {\it total Johnson homomorphism}
\index{total Johnson homomorphism of the mapping class group}
\index{Johnson homomorphism of the mapping class group!total}
of $\mathcal{M}_{g,1}$.
Then we have
\begin{theorem}[Hain \cite{Hai}]\label{T-Hain}
The Lie subalgebra $\mathrm{Im}(\tau_{\Q}^{\mathcal{M}})$ is generated by the degree one part
$\mathrm{Im}(\tau_{1,\Q}^{\mathcal{M}}) = \Lambda^3 H_{\Q}$ as a Lie algebra.
\end{theorem}

\vspace{0.5em}

Now, we consider the lower central series of the Torelli group, and reformulate Hain's result above.
For the lower central series $\{ \Gamma_{\mathcal{I}_{g,1}}(k) \}$ of $\mathcal{I}_{g,1}$, consider its graded quotients
$\mathrm{gr}^k(\mathcal{L}_{\mathcal{I}_{g,1}}) := \Gamma_{\mathcal{I}_{g,1}}(k)/\Gamma_{\mathcal{I}_{g,1}}(k+1)$ for $k \geq 1$.
Let
\[ {\tau'_{k}}^{\mathcal{M}} : \mathrm{gr}^k(\mathcal{L}_{\mathcal{I}_{g,1}}) \rightarrow H \otimes_{\Z} \mathrm{gr}^{k+1}(\mathcal{L}_{2g}) \]
be an $\mathrm{Sp}$-equivariant homomorphism induced from the restriction of $\tilde{\tau}_k$ to $\Gamma_{\mathcal{I}_{g,1}}(k)$. Then,
from Theorem {\rmfamily \ref{T-Hain}}, we have
\begin{proposition}[Hain, \cite{Hai}]\label{prop:Im}
$\mathrm{Im}(\tau_{k,\Q}^{\mathcal{M}})=\mathrm{Im}({\tau'_{k,\Q}}^{\hspace{-3.5mm}\mathcal{M}})$ for each $k \geq 1$.
\end{proposition}

Combining Hain's result above and the fact that $\mathrm{Coker}(\tau_{k,\Q}') \cong \mathcal{C}_n^{\Q}(k)$ for $n \geq k+2$ by Theorem {\rmfamily \ref{T-S11}},
we can detect non-trivial Sp-irreducible components in $\mathrm{Coker}(\tau_k^{\mathcal{M}})$.
In \cite{ES2}, we showed
\begin{theorem}[Enomoto and Satoh, \cite{ES2}]\label{ES-main}
For any $k \geq 5$ such that $k \equiv 1$ mod $4$, and $g \geq k+2$,
the irreducible Sp-module $[1^k]$ appears in $\mathrm{Coker}(\tau_{k,\Q}^{\mathcal{M}})$
with multiplicity one.
\end{theorem}

We wrote down a maximal vector of $[1^k]$ in \cite{ES2}.
In \cite{NT1}, Nakamura and Tsunogai gave a table of the Sp-irreducible decompositions of $\mathfrak{h}_{g,1}^{\Q}(k)$ for $1 \leq k \leq 15$,
which are calculated by a computer.
In the table, we can check that Sp-irreducible
components $[1^k]$ have multiplicity one for $k = 5, 9, 13$ and $k = 6, 10, 14$.
In 2004, Nakamura communicated to the author
the possibility that the multiplicities of $[1^{4m+1}]$ in $\mathfrak{h}_{g,1}(4m+1)$ is
exactly one for $m \geq 1$, and that they survive in the Johnson cokernels.
The theorem above is the affirmative answer to his conjecture.

\section{Twisted cohomology groups}\label{S-TwC}

In this section, we consider some twisted first cohomology groups of $\mathrm{Aut}\,F_n$ and the automorphism group of a free nilpotent group,
which are closely related to the study of the Johnson homomorphisms and trace maps.

\subsection{$H^1(\mathrm{Aut}\,F_n, \mathrm{IA}_n^{\mathrm{ab}})$}

First, we consider a twisted first cohomology groups of $\mathrm{Aut}\,F_n$ with coefficients in
$V:=\mathrm{IA}_n^{\mathrm{ab}} = H^* \otimes_{\Z} \Lambda^2 H$.
In \cite{S10}, we show that it is generated by two crossed homomorphisms
constructed with the Magnus representation and the Magnus expansion due to Morita and Kawazumi respectively.
As a corollary, we see that the first Johnson homomorphism does not extend to $\mathrm{Aut}\,F_n$ as a crossed homomorphism
for $n \geq 5$.

\vspace{0.5em}

Set
\[ \bm{e}_{j,k}^i := e_i^* \otimes e_j \wedge e_k \in V \]
for any $1 \leq i, j, k \leq n$, and fix them as a basis of $V$.
In \cite{S10}, using Nielsen's finite presentation for $\mathrm{Aut}\,F_n$, we computed
\begin{theorem}[Satoh, \cite{S10}]
For $n \geq 5$,
\[ H^1(\mathrm{Aut}\,F_n, V) = {\Z}^{\oplus 2}. \] 
\end{theorem}
This computation is a free group analogue of Morita's work \cite{Mo4} for the mapping class group $\mathcal{M}_{g,1}$.
In \cite{Mo4}, Morita computed the first cohomology group of $\mathcal{M}_{g,1}$ with coefficients in $\Lambda^3 H$, the free part of
the abelianization of the Torelli group $\mathcal{I}_{g,1}$, and showed
\begin{theorem}[Morita, \cite{Mo4}]
For $g \geq 3$,
\[ H^1(\mathcal{M}_{g,1},\Lambda^3 H)={\Z}^{\oplus 2}. \]
\end{theorem}

We give a description of a generators of $H^1(\mathrm{Aut}\,F_n, V)$.
To begin with, we construct a crossed homomorphism from $\mathrm{Aut}\,F_n$ into $V$ using the Magnus representation due to
Morita \cite{Mo1}.
We consider a left action of $\mathrm{Aut}\,F_n$.
Namely, define
\[ \bar{r}^{\mathfrak{a}} : \mathrm{Aut}\,F_n \longrightarrow \mathrm{GL}(n,\Z[H]) \]
by $\bar{r}^{\mathfrak{a}}(\sigma) := r^{\mathfrak{a}}(\sigma^{-1})$ for any $\sigma \in \mathrm{Aut}\,F_n$.
Then we have
\[ \bar{r}^{\mathfrak{a}}(\sigma \tau) = \bar{r}^{\mathfrak{a}}(\sigma) \cdot ( \sigma \cdot \bar{r}^{\mathfrak{a}}(\tau)) \]
for any $\sigma$, $\tau \in \mathrm{Aut}\,F_n$,
and hence $\bar{r}^{\mathfrak{a}}$ is also a crossed homomorphism (with respect to the left action).

\vspace{0.5em}

Let us recall Nielsen's finite presentation for $\mathrm{Aut}\,F_n$.
Let $P$, $Q$, $S$ and $U$ be automorphisms of $F_n$ given by specifying their images of the basis $x_1, \ldots, x_n$ as follows:
\begin{center}
\begin{tabular}{|c|c|c|c|c|c|c|} \hline
           & $x_1$      & $x_2$ & $x_3$ & $\cdots$ & $x_{n-1}$ & $x_n$ \\ \hline
  $P$      & $x_2$      & $x_1$ & $x_3$ & $\cdots$ & $x_{n-1}$ & $x_n$ \\ 
  $Q$      & $x_2$      & $x_3$ & $x_4$ & $\cdots$ & $x_{n}$   & $x_1$ \\ 
  $S$      & $x_1^{-1}$ & $x_2$ & $x_3$ & $\cdots$ & $x_{n-1}$ & $x_n$ \\ 
  $U$      & $x_1 x_2$  & $x_2$ & $x_3$ & $\cdots$ & $x_{n-1}$ & $x_n$ \\ \hline
\end{tabular}
\end{center}
In 1924, Nielsen \cite{Ni1} showed that the four elements above generate $\mathrm{Aut}\,F_{n}$.
Furthermore, he obtained the first finite presentation for $\mathrm{Aut}\,F_{n}$. (For details, see also \cite{MKS}.)

\vspace{0.5em}

By observing the images of Nielsen's generators by $\mathrm{det} \circ \bar{r}^{\mathfrak{a}}$, we verify that
$\mathrm{Im}(\mathrm{det} \circ \bar{r}^{\mathfrak{a}})$ is contained in
a multiplicative abelian subgroup $\pm H$ of $\Z[H]$. In order to modify the image of
$\mathrm{det} \circ \bar{r}^{\mathfrak{a}}$, we consider the signature of $\mathrm{Aut}\,F_n$.
For any $\sigma \in \mathrm{Aut}\,F_n$, set $\mathrm{sgn}(\sigma) := \mathrm{det}(\rho(\sigma)) \in \{ \pm 1 \}$,
and define a map $f_M : \mathrm{Aut}\,F_n \longrightarrow \Z[H]$
by
\[ \sigma \mapsto \mathrm{sgn}(\sigma) \,\, \mathrm{det}(\bar{r}^{\mathfrak{a}}(\sigma)). \]
Then the map $f_M$ is also a crossed homomorphism whose image is contained in a multiplicative abelian subgroup $H$ in $\Z[H]$.
In the following, we identify the multiplicative abelian group structure of $H$ with the additive one.

\vspace{0.5em}

In \cite{S01}, we computed a twisted first cohomology group of $\mathrm{Aut}\,F_n$ with coefficients in $H$ using Gersten's presentation
for $\mathrm{Aut}\,F_n$. That is,
\begin{theorem}[Satoh, \cite{S01}]
For $n \geq 2$,
\[ H^1(\mathrm{Aut}\,F_n, H) = \Z. \]
\end{theorem}
We showed that the crossed homomorphism $f_M$ generates $H^1(\mathrm{Aut}\,F_n, H)$.
This result is also a free group analogue of Morita's work in \cite{Mo0} for the mapping class group. That is,
\begin{theorem}[Morita, \cite{Mo0}]
For $g \geq 2$,
\[ H^1(\mathcal{M}_{g,1}, H) = \Z. \]
Furthermore, the restriction of the crossed homomorphism $f_M$ to $\mathcal{M}_{g,1}$ generates $H^1(\mathcal{M}_{g,1}, H)$.
\end{theorem}
This shows that the Dehn-Nielsen embedding $\mathcal{M}_{g,1} \hookrightarrow \mathrm{Aut}\,F_{2g}$ induces an isomorphism
\[ H^1(\mathrm{Aut}\,F_{2g},H) \xrightarrow{\cong} H^1(\mathcal{M}_{g,1},H). \]
See also Section 6 in \cite{KW2} by Kawazumi for details of a relation between the above twisted cohomology classes of $\mathcal{M}_{g,1}$ and the
Morita-Mumford classes.

\vspace{0.5em}

Now, let us turn back to the discussion of a generator of $H^1(\mathrm{Aut}\,F_n, \mathrm{IA}_n^{\mathrm{ab}})$.
Consider the inner automorphism group $\mathrm{Inn}\,F_n$ of $F_n$. Since $\mathrm{Inn}\,F_n$ is canonically isomorphic to $F_n$, its abelianization
can be identified with $H$.
By composing a homomorphism
$H \rightarrow V$ induced from the inclusion $\mathrm{Inn}\,F_n \hookrightarrow \mathrm{IA}_n$ and
the crossed homomorphism $f_M$, we obtain an element in $\mathrm{Cros}(\mathrm{Aut}\,F_n, V)$.
We also denote it by $f_M$.

\vspace{0.5em}

Next, we construct another crossed homomorphism from $\mathrm{Aut}\,F_n$ into $V$ using the Magnus expansion of $F_n$ by work of
Kawazumi \cite{Kaw}. (For basic material for the Magnus expansion, see Chapter 2 in \cite{Bou}, and also Section 7 in \cite{KW2}.)

\vspace{0.5em}

Let $\widehat{T}$ be the complete tensor algebra generated by $H$. For any Magnus expansion $\theta : F_n \rightarrow \widehat{T}$,
Kawazumi define a map
\[ \tau_1^{\theta} : \mathrm{Aut}\,F_n \rightarrow H^* \otimes_{\Z} H^{\otimes 2} \]
called the first Johnson map induced by the Magnus expansion $\theta$.
The map $\tau_1^{\theta}$ satisfies
\[ \tau_1^{\theta}(\sigma)([x]) = \theta_2(x) - |\sigma|^{\otimes 2} \theta_2(\sigma^{-1}(x)) \]
for any $x \in F_n$, where $[x]$ denotes the coset class of $x$ in $H$, $\theta_2(x)$ is the projection of $\theta(x)$ in $H^{\otimes 2}$,
and $|\sigma|^{\otimes 2}$ denotes the automorphism of $H^{\otimes 2}$ induced by $\sigma \in \mathrm{Aut}\,F_n$.
This shows that $\tau_1^{\theta}$ is a crossed homomorphism from $\mathrm{Aut}\,F_n$ to $H^* \otimes_{\Z} H^{\otimes 2}$.
In \cite{Kaw}, Kawazumi also showed that $\tau_1^{\theta}$ does not depend on the choice of the Magnus expansion $\theta$, and that
the restriction of $\tau_1^{\theta}$ to $\mathrm{IA}_n$ is a homomorphism satisfying
\[\begin{split}
  \tau_1^{\theta}(K_{ij}) & = e_i^* \otimes e_i \otimes e_j - e_i^* \otimes e_j \otimes e_i, \\
  \tau_1^{\theta}(K_{ijk}) & = e_i^* \otimes e_j \otimes e_k - e_i^* \otimes e_k \otimes e_j.
  \end{split}\]

\vspace{0.5em}

Now, compose a natural projection
$H^* \otimes_{\Z} H^{\otimes 2} \rightarrow H^* \otimes_{\Z} \Lambda^2 H =V$ and $\tau_1^{\theta}$.
Then we obtain an element in $\mathrm{Cros}(\mathrm{Aut}\,F_n,V)$. We denote it by $f_K$.
From the result of Kawazumi mentioned above, we see that the restriction of $f_K$ to $\mathrm{IA}_n$ coincides with
the double of the first Johnson homomorphism $\tau_1$. Namely, we have
\[ f_K(K_{ij}) = 2 \bm{e}_{i,j}^i = 2\tau_1(K_{ij}), \hspace{1em} f_K(K_{ijk}) = 2 \bm{e}_{j,k}^i=2\tau_1(K_{ijk}). \]

\vspace{0.5em}

Here, we consider the images of the crossed homomorphisms $f_M$ and $f_K$.
From the definition, we have
\[ f_M(\sigma) := \begin{cases}
                    -(\bm{e}_{1,2}^2 + \bm{e}_{1,3}^3 + \cdots + \bm{e}_{1,n}^n), \hspace{1em} & \sigma=S, \\
                    0, \hspace{1em} & \sigma=P, Q, U
                 \end{cases}\]
and
\[ f_K(\sigma) :=\begin{cases}
                    -\bm{e}_{1,2}^1, \hspace{1em} & \sigma=U, \\
                    0, \hspace{1em} & \sigma=P, Q, S.
                 \end{cases}\]
Moreover, in \cite{S10}, we showed that $f_M$ and $f_K$ generate $H^1(\mathrm{Aut}\,F_n, V_L)$ for $n \geq 5$. This shows that
the first Johnson homomorphism $\tau_1$ does not extend to a crossed homomorphism on $\mathrm{Aut}\,F_n$ for $n \geq 5$.

\vspace{0.5em}

It has already been known by Morita \cite{Mo5} that the first Johnson homomorphism
\[ \tau_1 : \mathcal{I}_{g,1} \rightarrow \Lambda^3 H \]
of the mapping class group does not extend to $\mathcal{M}_{g,1}$ as a crossed homomorphism for $g \geq 3$.

\subsection{$H^2(\mathrm{Aut}\,N_{n,k}, \Lambda^k H_{\Q})$ and a trace map $\mathrm{Tr}_{[1^k]}$}

Set
\[ \mathrm{Tr}_{[1^k]} := f_{[1^k]} \circ \Phi_{1}^k : H^* {\otimes}_{\Z} \mathcal{L}_n(k+1) \rightarrow {\Lambda}^k H. \]
We call it the {\it trace map for the exterior product} \index{trace map for the exterior product} $\Lambda^k H$.
In \cite{S03}, we show that $\Lambda^k H_{\Q}$ appears in $\mathrm{Coker}(\tau_{k,\Q}')$ for odd $k$ and $3 \leq k \leq n$, and
we determine $\mathrm{Coker}(\tau_{3,\Q})$ using $\mathrm{Tr}_{[3]}$ and $\mathrm{Tr}_{[1^3]}$.

\vspace{0.5em}

In \cite{ES1}, we prove that the trace map $\mathrm{Tr}_{[1^k]}$ defines a non-trivial twisted second cohomology class
of the automorphism group $\mathrm{Aut}\,N_{n,k}$ of a free nilpotent group $N_{n,k}:=F_n/\Gamma_n(k+1)$ with coefficients in
$\Lambda^k H_{\Q}$ for any $k \geq 2$ and $n \geq k$. Here we give a brief strategy to prove it.

\vspace{0.5em}

First, we recall Bryant and Gupta's generators of $\mathrm{Aut}\,N_{n,k}$.
Let $\mathrm{Aut}\,F_n \rightarrow \mathrm{Aut}\,N_{n,k}$ be a homomorpshim induced from the natural projection $F_n \rightarrow N_{n,k}$.
Andreadakis \cite{And} showed that it is surjective if $k=2$, and not if otherwise.
For any $\sigma \in \mathrm{Aut}\,F_n$, we also denote its image in $\mathrm{Aut}\,N_{n,k}$ by $\sigma$.
Under this notation, $\mathrm{Aut}\,N_{n,2}$ is generated by Nielsen's generators $P$, $Q$, $S$ and $U$.

\vspace{0.5em}

Let $\theta$ be an automorphism of $N_{n,k}$, defined by 
\[ \theta : x_t \mapsto \begin{cases}
               [x_1, [x_2,x_1]]x_1, & t=1, \\
               x_t,                & t \neq 1.
              \end{cases}\]
Then we have
\begin{theorem}[Bryant and Gupta \cite{BrG}]
For $k \geq 3$ and $n \geq k-1$, the group $\mathrm{Aut}\,N_{n,k}$ is generated by
$P$, $Q$, $S$, $U$ and $\theta$.
\end{theorem}
We remark that in 1984, Andreadakis \cite{An2} showed that $\mathrm{Aut}\,N_{n,k}$ is generated by $P$, $Q$, $S$, $U$ and $k-2$ other elements for $n \geq k \geq 2$.
No presentation for $\mathrm{Aut}\,N_{n,k}$ is known except for $\mathrm{Aut}\,N_{2,k}$ for $k= 1, 2$ and $3$
due to Lin \cite{Lin}.

\vspace{0.5em}

Using Bryant and Gupta's generators and some relations among them, we showed
\begin{proposition}[Enomoto and Satoh, \cite{ES1}]
For $k \geq 3$, $n \geq k-1$ and $l \geq 3$,
$H^1(\mathrm{Aut}\,N_{n,k}, \Lambda^l H_{\Q})$ is trivial.
\end{proposition}
Then observing a group extension
\[ 0 \rightarrow \mathrm{Hom}_{\Z}(H, \mathcal{L}_n(k+1)) \rightarrow \mathrm{Aut}\, N_{n,k+1} \rightarrow \mathrm{Aut}\,N_{n,k} \rightarrow 1 \]
introduced by Andreadakis \cite{And}, (See also Proposition 2.3 in Morita's paper \cite{Mo4}.),
and its cohomological five term exact sequence, we obtain
\[\begin{split}
   0 & \rightarrow H^1(H^* \otimes_{\Z} \mathcal{L}_n(k+1), \Lambda^k H_{\Q})^{\mathrm{GL}(n,\Z)} \\
     & \rightarrow H^2(\mathrm{Aut}\,N_{n,k}, \Lambda^k H_{\Q}) \rightarrow H^2(\mathrm{Aut}\,N_{n,k+1}, \Lambda^k H_{\Q})
  \end{split}\]
for $k \geq 3$ and $n \geq k-1$. 

\vspace{0.5em}

Since the trace map
$\mathrm{Tr}_{[1^k]} \in H^1(H^* \otimes_{\Z} \mathcal{L}_n(k+1), \Lambda^k H_{\Q})^{\mathrm{GL}(n,\Z)}$
is surjective for any $3 \leq k \leq n$, we have
\begin{proposition}
For $k \geq 3$ and $n \geq k$, we see
\[ 0 \neq \mathrm{tg}(\mathrm{Tr}_{[1^k]}) \in H^2(\mathrm{Aut}\,N_{n,k}, \Lambda^k H_{\Q}) \]
where $\mathrm{tg}$ is the transgression map.
\end{proposition}

\vspace{0.5em}

In \cite{ES1}, we considered not only $\mathrm{Aut}\,N_{n,k}$ but also
the image of the natural homomorphism $\mathrm{Aut}\,F_n \rightarrow \mathrm{Aut}\,N_{n,k}$, denoted by $T_{n,k}$.
The group $T_{n,k}$ is called the {\it tame automorphism group} \index{tame automorphism group} of $N_{n,k}$. Similarly to $\mathrm{Aut}\,N_{n,k}$,
observing the cohomological five term exact sequence of a group extension
\[ 0 \rightarrow \mathrm{gr}^k(\mathcal{A}_n) \rightarrow T_{n,k+1} \rightarrow T_{n,k} \rightarrow 1, \]
we showed
\begin{proposition}
For even $k$ and $2 \leq k \leq n$, we see
\[ 0 \neq \mathrm{tg}(\mathrm{Tr}_{[1^k]} \circ \tau_{k}) \in H^2(T_{n,k}, \Lambda^k H_{\Q}) \]
where $\mathrm{tg}$ is the transgression map.
\end{proposition}

\vspace{0.5em}

In \cite{ES1}, we showed that
\[ \mathrm{Hom}_{\Z}(H_{\Q}^* \otimes_{\Z} \mathcal{L}_{n,\Q}(k+1), \Lambda^k H_{\Q})^{\mathrm{GL}(n,\Q)}
    \cong \Q \]
for $n \geq 2k+2$, and that
\[ \mathrm{Hom}_{\Z}(\mathrm{gr}_{\Q}^k(\mathcal{A}_n), \Lambda^k H_{\Q})^{\mathrm{GL}(n,\Q)} \cong \Q \]
for even $k$ and $n \geq 2k+2$.
At the present stage, however, it seems too difficult to determine the precise structures of
$H^1(H^* \otimes_{\Z} \mathcal{L}_n(k+1), \Lambda^k H_{\Q})^{\mathrm{GL}(n,\Z)}$ and
$H^1(\mathrm{gr}^k(\mathcal{A}_n), \Lambda^k H_{\Q})^{\mathrm{GL}(n,\Z)}$ in general.

\section{Congruence IA-automorphism groups of $F_n$}\label{S-Cong}

For $n \geq 2$ and $d \geq 2$, let
$\mathrm{GL}(n,\Z) \rightarrow \mathrm{GL}(n,\Z/d\Z)$ be the natural homomorphism induced by the mod reduction $d$.
The kernel of this homomorphism is called the congruence subgroup of $\mathrm{GL}(n,\Z)$ of level $d$, and is denoted by $\Gamma(n,d)$.
Let $\mathrm{IA}_{n,d}$ be the kernel of the composite homomorphism
$\mathrm{Aut}\,F_n \xrightarrow{\rho} \mathrm{GL}(n,\Z) \rightarrow \mathrm{GL}(n,\Z/d\Z)$.
We call $\mathrm{IA}_{n,d}$ the {\it congruence IA-automorphism group} \index{congruence IA-automorphism group}
of $F_n$ of level $d$. Then we have a group extension
\begin{equation}\label{ex-cia}
 1 \rightarrow \mathrm{IA}_n \rightarrow \mathrm{IA}_{n,d} \rightarrow \Gamma(n,d) \rightarrow 1.
\end{equation}
In this section, we give a description of the abelianization of $\mathrm{IA}_{n,d}$ with those of $\mathrm{IA}_n$ and $\Gamma(n,d)$
by using the first Johnson map.

\subsection{The abelianization of $\mathrm{IA}_{n,d}$}
\vspace{0.5em}

In \cite{S02} and \cite{S02'}, using the first Johnson map over $\Z/d\Z$ defined by the Magnus expansion due to Kawazumi \cite{Kaw}, we showed
\begin{theorem}[Satoh, \cite{S02} and \cite{S02'}]\label{T-KaS}
For $n \geq 2$ and $d \geq 2$,
\[ \mathrm{IA}_{n,d}^{\mathrm{ab}} \cong (\mathrm{IA}_n^{\mathrm{ab}} \otimes_{\Z} \Z/d\Z) \, \bigoplus \, {\Gamma(n,d)}^{\mathrm{ab}}. \]
\end{theorem}

We give the outline of the proof. For simplicity, we consider the case where $d$ is odd.
Observing the homological five term exact sequence of the group extension (\ref{ex-cia}), we obtain an exact sequence
\[\begin{split}
   \cdots & \rightarrow H_1(\mathrm{IA}_n, \Z)_{\Gamma(n,d)} \xrightarrow{i} H_1(\mathrm{IA}_{n,d},\Z) 
     \rightarrow H_1(\Gamma(n,d),\Z) \rightarrow 0.
  \end{split}\]
Then we show that
\[ H_1(\mathrm{IA}_n, \Z)_{\Gamma(n,d)} \cong H_1(\mathrm{IA}_n,\Z) \otimes_{\Z} \Z/d\Z, \]
and that the first Johnson map
\[ \tau^{\theta} : H_1(\mathrm{IA}_{n,d},\Z) \rightarrow H_1(\mathrm{IA}_n,\Z) \otimes_{\Z} \Z/d\Z \]
satisfies $i \circ \tau^{\theta} = \mathrm{id}$. This means that
\[ 0 \rightarrow H_1(\mathrm{IA}_n, \Z)_{\Gamma(n,d)} \xrightarrow{i} H_1(\mathrm{IA}_{n,d},\Z) 
     \rightarrow H_1(\Gamma(n,d),\Z) \rightarrow 0 \]
is a split exact sequence. This shows Theorem {\rmfamily \ref{T-KaS}}.
If $d$ is even, we need more complicated arguments since the image of $\tau^{\theta}$ is not contained in
$H_1(\mathrm{IA}_n,\Z) \otimes_{\Z} \Z/d\Z$. (See \cite{S02'} for details.)

\vspace{0.5em}

If $d$ is odd prime $p$, then the structure of $\Gamma(n,p)^{\mathrm{ab}}$ is well-known
due to Lee and Szczarba.
Set
\[ \mathfrak{sl}(n, \Z/p\Z) := \{ A \in M(n, \Z/p\Z) \,|\, \mathrm{Trace}\,A = 0 \}. \]
Then we have
\begin{theorem}[Lee and Szczarba \cite{LeS}]
For $n \geq 3$ and odd prime $p$, $\Gamma(n,p)^{\mathrm{ab}} \cong \mathfrak{sl}(n,\Z/p\Z)$ as an $\mathrm{SL}(n,\Z/p\Z)$-module.
\end{theorem}

From this theorem, we see that for any odd prime $p$, the abelianization
of $\mathrm{IA}_{n,p}$ is isomorphic to $(\Z/p\Z)^{\oplus \frac{1}{2}(n-1)(n^2+2n+2)}$ as an abelian group.

\vspace{0.5em}

We remark that a mapping class group analogue of the above was shown independently by Putman \cite{Put} and Sato \cite{MaS}.
Let $\mathcal{I}_{g,1}(d)$ be the kernel of the composite homomorphism
$\mathcal{M}_{g,1} \rightarrow \mathrm{Sp}(2g,\Z) \rightarrow \mathrm{Sp}(2g,\Z/d\Z)$.
Set
\[ \mathfrak{sp}(2g, \Z/d\Z) := \{ A \in M(2g, \Z/d\Z) \,|\, {}^t A J + JA = O \}. \]
Putman and Sato showed
\begin{theorem}[Putman \cite{Put} and \cite{Pu2}, and Sato \cite{MaS}]
For any $g \geq 3$ and odd $d \geq 3$,
\[ \mathcal{I}_{g,1}(d)^{\mathrm{ab}} \cong (\mathcal{I}_{g,1}^{\mathrm{ab}} \otimes_{\Z} \Z/d\Z) \, \bigoplus \, \mathfrak{sp}(2g, \Z/d\Z). \]
\end{theorem}
This result is generalized by Putman \cite{Pu2} for an even $d \geq 2$ which is not divisible by $4$.
Sato \cite{MaS} also proved the above for $d=2$ independently. That is,
\begin{theorem}[Putman \cite{Pu2} and Sato \cite{MaS}]
For any $g \geq 3$,
\[ \mathcal{I}_{g,1}(2)^{\mathrm{ab}} \cong (\Z/2\Z)^{\oplus \binom{2g}{3}} \bigoplus (\Z/4\Z)^{\oplus \binom{2g}{2}} \bigoplus
   (\Z/8\Z)^{\oplus \binom{2g}{1}}. \]
\end{theorem}

\subsection{The case of $n=2$}

For $n=2$ and any odd prime $p$, the group structure of $\mathrm{IA}_{2,p}$ is much easier to handle than those for $n \geq 3$.
We have a group extension
\begin{equation}\label{ex-cia2}
 1 \rightarrow \mathrm{Inn}\,F_2 \rightarrow \mathrm{IA}_{2,p} \rightarrow \Gamma(2,p) \rightarrow 1.
\end{equation}
By a classical result by Frasch, we have
\begin{theorem}[Frasch, \cite{Fra}]
For any odd prime $p$, the congruence subgroup $\Gamma(2,p)$ is a free group of rank $\alpha(p):= 1+ p(p^2-1)/12$.
\end{theorem}
Hence, $\mathrm{IA}_{2,p}$ is a group extension of free groups of finite rank.

\vspace{0.5em}

In \cite{S02}, we compute the integral homology groups of $\mathrm{IA}_{2,p}$ for an odd prime $p$, using the
Lyndon-Hochshild-Serre spectral sequence
of the extension (\ref{ex-cia2}).
\begin{proposition}[Satoh, \cite{S02}]
For any prime $p$,
\[ H_q(\mathrm{IA}_{2,p},\Z) =
              \begin{cases}
               \Z & \quad \mathrm{if} \,\,\, q=0, \\
               {\Z}^{\oplus \alpha(p)} \oplus (\Z/p\Z)^{\oplus 2} & \quad \mathrm{if} \,\,\, q=1,\\
               {\Z}^{\oplus (2\alpha(p)-2)} & \quad \mathrm{if} \,\,\, q=2,\\
               0 &\quad\mathrm{if} \,\,\, q \geq 3.
              \end{cases} \]
\end{proposition}

\section{Acknowledgments}

The author would like to thank Professor Athanase Papadopoulos and Professor Nariya Kawazumi
for giving me an opportunity to write this chapter, and to contribute to the Handbook of Teichm$\ddot{\mathrm{u}}$ller Theory.
He thanks Professor Athanase Papadopoulos for careful reading and valuable comments on this chapter.
He also thanks to Andrew Putman for communicating his recent works with respect to the abelianization of $\mathcal{I}_{g,1}(d)$.

\end{document}